\documentclass{article}

\usepackage[utf8]{inputenc}
\usepackage[english]{babel}
\usepackage{amsfonts}
\usepackage[11pt]{moresize}
\usepackage{pgfplots}
\usepackage{xfrac}
\usepackage{lipsum}

\newcommand\freefootnote[1]{%
	\let\thefootnote\relax%
	\footnotetext{#1}%
	\let\thefootnote\svthefootnote%
}

\newcommand{\Addresses}{{
		\bigskip
		\footnotesize
		
		Humboldt-Universit\"at zu Berlin, Institut f\"ur Mathematik, Rudower Chausee 25
			\hfill \newline\texttt{}
			\indent 12489 Berlin, Germany} 
		\par\nopagebreak
		\textit{E-mail address}: \texttt{andreibud95@protonmail.com}
	}

\usepackage{amssymb}
\usepackage{mathtools}
\usepackage{tikz}
\usepackage{tikz-cd}
\usepackage{mathtools}
\usepackage{amsmath}
\usepackage{amsthm}
\usepackage{hyperref}
\usepackage{float}
\usetikzlibrary{calc, positioning,decorations.markings,arrows}

\usepackage{graphics}
\usepackage{color}
\usepackage{tabularx,colortbl}
\usetikzlibrary{positioning}
\usetikzlibrary{decorations.markings,arrows}
\usetikzlibrary{calc}
\usepackage{amsthm}

\theoremstyle{plain}
\newtheorem{trm}{Theorem}[section]
\newtheorem{lemma}[trm]{Lemma}
\newtheorem{proposition}[trm]{Proposition}
\newtheorem{corollary}[trm]{Corollary}
\newtheorem{theorem}[trm]{Theorem}

\theoremstyle{definition}

\newtheorem{remark}[trm]{Remark}

\def\Pic0{{\rm Pic}^0(X)}

\begin{document}
	
	\title{A Hurwitz divisor on the moduli of Prym curves}

	
	\author{Andrei Bud
	}
\date{}
	
	\maketitle
	
	\begin{abstract}
		For even genus $g=2i\geq4$ and the length $g-1$ partition $\mu = (4,2,\ldots,2,-2,\ldots,-2)$ of 0, we compute the first coefficients of the class of $\overline{D}(\mu)$ in $\mathrm{Pic}_\mathbb{Q}(\overline{\mathcal{R}}_g)$, where $D(\mu)$ is the divisor consisting of pairs $[C,\eta]\in \mathcal{R}_g$ with $\eta \cong \mathcal{O}_C(2x_1+x_2+\cdots + x_{i-1}-x_i-\cdots-x_{2i-1})$ for some points $x_1,\ldots, x_{2i-1}$ on $C$. We further provide several enumerative results that will be used for this computation. \par \  
		\end{abstract}
	
	\section{Introduction}  \label{sec1} \freefootnote{\Addresses}
	The moduli space $\mathcal{R}_g$ parametrizing pairs $[C,\eta]$ consisting of a curve $C$ of genus $g$ and a $2$-torsion line bundle $\eta$ on $C$ received considerable attention following the influential papers \cite{Beau77} and \cite{MumfordPrym}. The description of $\mathcal{R}_g$ in \cite{Beau77} as a coarse moduli space of a stack, together with the algebraic theory of Prym curves developed by Mumford in \cite{MumfordPrym} brought this topic to the attention of algebraic geometers. To outline its importance, we recall that $\mathcal{R}_g$ comes equipped with a map $\mathcal{P}_g\colon \mathcal{R}_g \rightarrow \mathcal{A}_{g-1}$ to the moduli space of principally polarized abelian varieties of dimension $g-1$. This natural application relating curves to Prym varieties inside $\mathcal{A}_{g-1}$ was used to provide an algebraic proof of the Schottky-Jung relations, see \cite{MumfordPrym}, and, among others, to understand the birational geometry of the moduli of Prym varieties, see \cite{FarLud}, \cite{Bruns}, \cite{FarVerNikulin} and the references therein. 
	
	Let $\overline{\mathcal{R}}_g$ be the compactification of $\mathcal{R}_g$ as considered in \cite{Beau77} and \cite{Casa}. When $g =2i+1$ is odd, Farkas and Ludwig considered an effective divisor $\overline{D}_{2i+1:2} = \overline{\left\{[C,\eta]\in\mathcal{R}_{2i+1} \ | \ \eta \in C_i-C_i \right\}} $ describing the relative position of $\eta$ with respect to the divisor $C_i-C_i$ in $\mathrm{Pic}^0(C)$, and computed some relevant coefficients of its class in $\mathrm{Pic}_\mathbb{Q}(\overline{\mathcal{R}}_g)$. As a consequence of their computation, they obtained that $\mathcal{R}_g$ is of general type for $g=13, 14$ and $g\geq 17$. A natural adaptation of their treatment to the case when $g=2i$ is even is to consider the divisor in $\mathcal{R}_{g}$ parametrizing pairs $[C,\eta]$ satisfying $\eta \cong \mathcal{O}_C(2x_1+x_2+\cdots+x_{i-1}-x_i-x_{i+1}-\cdots-x_{g-1})$ for some points $x_1,\ldots, x_{g-1}$ on $C$, and compute some of its coefficients in $\mathrm{Pic}_\mathbb{Q}(\overline{\mathcal{R}}_g)$. In terms of the position of $\eta$ with respect to a difference divisor in $\mathrm{Pic}^0(C)$, the divisor we study is $\overline{\left\{[C,\eta]\in \mathcal{R}_{2i} \ | \ \eta\in 2C+ C_{i-2}-C_i\right\}}$. 
	
	Looking at this from another perspective, we can see the divisor $D_{2i+1:2}$ as the image in $\mathcal{R}_{2i+1}$ of some Hurwitz scheme. While for the moduli space of curves $\mathcal{M}_g$ several cases when a Hurwitz locus is a divisor are studied, see \cite{KodMg}, \cite{KodevenHarris1984}, \cite{Diaz}, \cite{FarkasFermat} and \cite{VGeerKouvid}, on the moduli space $\mathcal{R}_g$ the only studied example is $D_{2i+1:2}$. Such Hurwitz divisors are fundamental in proving that $\mathcal{M}_g$ is of general type for $g\geq 24$, see \cite{KodMg}, \cite{KodevenHarris1984} and \cite{EisenbudHarrisg>23}, and that $\mathcal{R}_g$ is of general type for $g\geq17$, see \cite{FarLud}.    

	
	Let $\mu = (2m_1,2m_2,\ldots, 2m_{g-1})$ a length $g-1$ partition of $0$. We can define a Hurwitz divisor, denoted $D(\mu)$ as the locus parametrizing pairs $[C, \eta]$ satisfying a line bundle isomorphism $\eta \cong \mathcal{O}_C(m_1x_1+\cdots +m_{g-1}x_{g-1})$ for some points $x_1,\ldots, x_{g-1}$ of $C$. When $g=2i+1$ and $\mu$ is $(2,\ldots,2,-2,\ldots,-2)$ we recover the divisor $D_{2i+1:2}$. In this article we are interested in the case $g=2i$ with partition $\mu = (4,2,\ldots,2,-2,\ldots,-2)$ and we ask what is the class of the divisor $\overline{D}(\mu)$ in $\mathrm{Pic}_\mathbb{Q}(\overline{\mathcal{R}}_g)$. We consider the basis of $\mathrm{Pic}_\mathbb{Q}(\overline{\mathcal{R}}_g)$ consisting of the classes $\lambda, \delta_0', \delta_0'', \delta_0^{\mathrm{ram}}$ together with $\delta_i, \delta_{g-i}, \delta_{i:g-i}$ for $1 \leq i \leq [g/2]$. In this basis, we compute the first coefficients of our divisor $\overline{D}(\mu)$ and we obtain: 
	
	\begin{theorem}\label{maintheorem} Let $g =2i\geq4$ and $\mu$ the length $g-1$ partition $(4,2,\ldots,2,-2,\ldots,-2)$ of $0$. Then for the class in $\mathrm{Pic}_\mathbb{Q}(\overline{\mathcal{R}}_g)$ of the divisor 
		\[ [\overline{D}(\mu)]  \equiv  a\lambda - b_0'\delta_0' - b_0''\delta_0'' - b_0^{\mathrm{ram}}\delta_0^{\mathrm{ram}}- b_1\delta_1 - b_{g-1}\delta_{g-1} - b_{1:g-1}\delta_{1:g-1}- \cdots  \]
		we have the equalities: 
		\[ a = \frac{12i^2+10i-2}{2i-1}\cdot\binom{2i-1}{i}, \ \ b_0' = \frac{2i^2}{2i-1}\cdot\binom{2i-1}{i}\]
		\[ b_0'' = \frac{4i^3}{2i-1}\cdot\binom{2i-1}{i}-(3i-1)\cdot 2^{2i-2}, \ \ b_0^{\mathrm{ram}} = \frac{2i^2+3i-1}{2i-1}\cdot\binom{2i-1}{i}\]
		\[b_1 = 2i(4i+1)\cdot\binom{2i-1}{i}-6(2i-1)\cdot 2^{2i-2}, \ \ b_{g-1} = (6i-2)\binom{2i-1}{i}, \ \  b_{1:g-1} = (2i+2)\cdot\binom{2i-1}{i} \]
	\end{theorem}
	
	Due to results in \cite{FarkasPopaK3} and \cite{FarLud}, when $g\leq 23$ all the other coefficients are irrelevant from the point of view of birational geometry. As many interesting questions about the birational type appear when $g\leq 23$, no relevant information is lost in this way. 
	
	In order to prove this result, we will consider the compactification of Hurwitz schemes by means of admissible covers, see \cite{KodMg}, \cite{Diaz} and \cite{AbramovichCV}. Further, we intersect our divisor with some classical test curves and compute the number of admissible covers above this intersection, along with their multiplicities. We get in this way a system of $8$ equations with $7$ unknowns which will be compatible and will conclude Theorem \ref{maintheorem}.
	
	\section{Admissible covers and enumerative geometry} 
	We begin by providing the setting, along with some important results about admissible covers. Next, we present the enumerative results we need in order to compute the intersection of our divisor with different test curves. 
	
	Let $\mu = (2m_1, 2m_2,\ldots, 2m_{g-1})$ a length $g-1$ partition of 0 and let $\mu^{-}$ and $\mu^{+}$ be the vectors of negative and positive entries of $\mu$. We denote by $d$ the sum of the positive elements in $\mu$ and take the partitions of $d$ given as $b_1 = \mu^{+}$, $b_2= -\mu^{-}, b_3 = \cdots = b_{3g-1} = (2,1\ldots, 1)$. We further consider the set $B=\left\{b_1,b_2, \ldots b_{3g-1} \right\}$. 
	
	Following the notation in \cite{Diaz}, we consider the moduli space $H_{d,B}$ parametrizing degree $d$ maps $[\pi\colon X \rightarrow \mathbb{P}^1]$ together with points $q_1, q_2, \ldots, q_{3g-1}$ on $\mathbb{P}^1$ such that over $q_i$ the map has ramification profile $b_i$ and is otherwise unramified. As previously mentioned, we have a compactification $\overline{H}_{d,B}$ of $H_{d,B}$ by means of admissible covers. We remark that the complete local rings of the Hurwitz scheme $\overline{H}_{d,B}$ are determined as in \cite{KodMg}.
	
	Let $\mathcal{S}(\mu)$ be the subgroup of $S_{g-1}$ generated by the transpositions $\left\{(i,j) \ | \ m_i = m_j\right\}$. This group acts on $\overline{\mathcal{M}}_{g,g-1}$ by permuting the marked points, and we can consider the quotient $\overline{\mathcal{M}}_{g,g-1}/\mathcal{S}(\mu)$ of this action. We then have the map 
	\[ a_\mu\colon\overline{H}_{d,B} \rightarrow \overline{\mathcal{M}}_{g,g-1}/\mathcal{S}(\mu)\]
	sending $[\pi\colon X \rightarrow \mathbb{P}^1]$ with branch points $q_1, q_2,\ldots, q_{3g-1}$ to the stable model of the pointed curve $[X, p_1,\ldots,p_{g-1}]$ where $p_1,\ldots, p_{g-1}$ are the points in the preimages of $q_1$ and $q_2$, considered in the order given by $\mu$
	
	We consider the subspace of $\mathcal{M}_{g,g-1}$ defined as \[\mathcal{H}^0_g(\mu) = \left\{[X,p_1,\ldots, p_{g-1}] \in \mathcal{M}_{g,g-1}  \mid \mathcal{O}_X\left(\sum_{i=1}^{g-1}2m_ip_i\right) \cong \mathcal{O}_X \right\}\] 
	and we observe that the image of the map $a_\mu$ is $\overline{\mathcal{H}}_g^0(\mu)/\mathcal{S}(\mu)$.

	We want to discount the components $\mathcal{H}^0_g(\frac{\mu}{2})$ in the space $\mathcal{H}^0_g(\mu)$ and hence, we will only consider the components of $H_{d,B}$ mapping to $\text{\large (}\mathcal{H}^0_g(\mu) \setminus \mathcal{H}^0_g(\frac{\mu}{2})\text{\large )} /\mathcal{S}(\mu)$. We will denote by $H_{d,\mu}$ the space of these components.
	
	This restriction enables us to consider a map 
	\[ b_\mu\colon \left(\mathcal{H}^0_g(\mu) \setminus \mathcal{H}^0_g(\frac{\mu}{2}) \right)/\mathcal{S}(\mu) \rightarrow \mathcal{R}_g \]
	given as 
	\[ b_\mu\text{\large (} [X,p_1,\ldots,p_{g-1}] \text{\large )} = [X, \mathcal{O}_X(\sum_{i=1}^{g-1}m_ip_i)] \]
	
	Using the compactification of $\mathcal{H}^0_g(\mu)$ in terms of twists at the nodes, see \cite{FP18}, we are able to extend the map $b_\mu$ over curves of compact type: A pointed curve $[X,p_1,\ldots, p_{g-1}]$ is sent to $[\mathrm{stab}(X), \eta]$ where if $C$ is a component of the stabilization of $X$, we take $\eta_C$ to be $\mathcal{O}_C(\sum n_iq_i)$ where $\sum 2n_iq_i = 0$ is the divisorial equivalence determined by the unique possible twist on the component $C$. 
	
	Let $c_\mu = b_\mu\circ a_\mu$ and $D(\mu) = c_{\mu*}(H_{d,\mu})$. This is a divisor in $\mathcal{R}_g$ and we want to compute the class of its closure $\overline{D}(\mu)$ in $\overline{\mathcal{R}}_g$. Let $\pi_\mu\colon \overline{H}_{d,\mu} \rightarrow \overline{\mathcal{M}}_g$ be the map sending an admissible cover $[\pi\colon X\rightarrow \Gamma]$ to the curve $\mathrm{stab}(X) \in \overline{\mathcal{M}}_g$. We consider the divisor $\overline{Z}(\mu) = \pi_{\mu *}(\overline{H}_{d,\mu})$. Then it is obvious that the set-theoretical projection of $\overline{D}(\mu)$ to $\overline{\mathcal{M}}_g$ is $\overline{Z}(\mu)$. This observation implies the following:
	
	\begin{lemma} \label{lm1}
		Let $[C,\eta]$ a point in the divisor $\overline{D}(\mu)$. Then there exists an admissible cover $[\pi\colon X \rightarrow \Gamma]$ in the Hurwitz scheme $\overline{H}_{d,\mu}$ such that $\mathrm{StMd}(X) = \mathrm{StMd}(C)$.  
	\end{lemma} 
	
	Our goal is to particularize to $g=2i$ and $\mu = (4,2,\ldots, 2, -2,\ldots, -2)$, and prove Theorem \ref{maintheorem}. We make the following immediate observation which implies that we can ignore the components of $H_{d,B}\setminus H_{d,\mu}$ in our computations.
	\begin{remark} Let $\mu = (2m_1,\ldots, 2m_{g-1})$ a partition of 0 with all negative entries equal to $-2$. If $Z$ is a component of $H_{d,B}$ mapped to $\mathcal{H}^0_{g}(\frac{\mu}{2})/\mathcal{S}(\mu)$ by $a_\mu$, then the projection of $Z$ to $\mathcal{M}_g$ has at least codimension 2. 
	\end{remark}
	
	\subsection{Various enumerative results} 
	
	We provide here some enumerative results that will be used to compute the intersection of the divisor $\overline{D}(\mu)$ for $\mu = (4,2,\ldots,2,-2,\ldots,-2)$ with various test curves. We start by counting the number of maps to $\mathbb{P}^1$ satisfying some ramification conditions on a special fibre. Subsequently, we switch our attention to elliptic curves and compute the degrees of some particular Hurwitz schemes over $\mathcal{M}_{1,1}$.
	
	We state and prove a classical result generalizing Theorem B in \cite{KodMg}, Theorem 2.1 in \cite{KodevenHarris1984} and Lemma 6.2 in \cite{Diaz}. 
	
	\begin{theorem}\label{BrillNoether} For $m,n$ non-negative integers, let $d\leq g+1-m$, $[C,x_1,\ldots, x_n]$ a generic point in $\mathcal{M}_{g,n}$ and $\alpha_1,\ldots, \alpha_m, \beta_1,\ldots, \beta_n$ positive integers satisfying 
		\[ \sum_{i=1}^{m} \alpha_i + \sum_{j=1}^{n} \beta_j = 2d+m-1-g. \]
		Then the number of pairs $(L,y_1,\ldots,y_m)$ with $L$ a degree $d$ line bundle on $C$ and $y_1,\ldots, y_m$ points on $C$ satisfying 
		\[ h^0(C, L) \geq 2 \mathrm{\ and \ } h^0(C, L(-\sum_{i=1}^{m}\alpha_iy_i - \sum_{j=1}^{n}\beta_jx_j)) \geq 1\]
		is equal to $N_{\alpha,\beta} \coloneqq \frac{g!}{d!(g+1-d-m)!} (2d+m-g-1-\sum_{i=1}^{m}\frac{1}{\alpha_i})\prod_{i=1}^{m}\alpha_i^2$.
		
		Moreover, due to the genericity of $[C,x_1,\ldots, x_n]$, for every such pair $(L,y_1,\ldots, y_m)$ the line bundle $L$ is globally generated, $h^0(C, L) = 2$ and $h^0(C, L(-\sum_{i=1}^{m}\alpha_iy_i - \sum_{j=1}^{n}\beta_jx_j)) = 1$. Furthermore, the points $x_1,\ldots, x_n, y_1,\ldots, y_m$ along with the points in the support of the effective divisor $D$ equivalent to $\mathrm{div}(L)-\sum_{i=1}^{m}\alpha_iy_i - \sum_{j=1}^{n}\beta_jx_j$ are pairwise distinct.   
	\end{theorem}
	\begin{proof} First, observe that everything after "moreover" is true by dimension considerations.
		Consider the map 
		\[ \varphi_{\alpha,\beta}\colon \underbrace{C\times\cdots\times C}_{\textrm{ $m$ times}}\times C_{g+1-d-m} \rightarrow C_d  \]
		given as 
		\[ (y_1,\ldots, y_m, D) \mapsto D +\sum_{i=1}^m \alpha_iy_i+ \sum_{j=1}^{n}\beta_jx_j \]
		Let $C^1_d \subseteq C_d$ be the locus of divisors $D$ satisfying $h^0(C, D)\geq 2$. It is easy to observe that the number of pairs $(L,y_1,\ldots,y_m)$ satisfying the desired properties is the number of points in the intersection $\mathrm{Im}(\varphi_{\alpha,\beta})\cap C^1_d$.
		
		We will show that this intersection is transverse. Consider a point 
		\[ E = \sum_{i=1}^m\alpha_iy_i + \sum_{j=1}^n \beta_jx_j + \sum_{k=1}^{g+1-d-m}z_k \]
		in the intersection $\mathrm{Im}(\varphi_{\alpha,\beta})\cap C^1_d$. Because the points are all distinct and $h^0(C, E) = 2$ it follows that both $\mathrm{Im}(\varphi_{\alpha,\beta})$ and $C^1_d$ are smooth at the point $E$. 
		Inside the tangent space 
		\[T_E(C_d) = H^0(C, \mathcal{O}_C(E)/\mathcal{O}_C) = H^0(C,\omega_C/\omega_C(-E))^\vee\] 
		we have the following identifications
		\[ T_E(\mathrm{Im}(\varphi_{\alpha,\beta})) = \textrm{ Annihilator of\ } H^0\textrm{\Large (}C, \omega_C\textrm{\large (}-E + \sum_{i=1}^m (\alpha_i -1)y_i + \sum_{j=1}^{n}\beta_jx_j\textrm{\large )}/\omega_C(-E)\textrm{\Large )} \]
		\[ T_E(C^1_d) = \textrm{ Annihilator of \ } \mathrm{Im}(\tilde{\mu}_0)  \]
		where $\tilde{\mu}_0$ is the composition 
		\[ H^0(C,\mathcal{O}_C(E))\otimes H^0(C, \omega_C(-E))\rightarrow H^0(C,\omega_C) \rightarrow H^0\textrm{\large(}C, \omega_C/\omega_C(-E)\textrm{\large)} \]
		It is clear that for any differential $s \in H^0(C,\omega_C(-E))$ we have $\tilde{\mu}_0(1\otimes s) = 0$. Consider $f$ in $H^0(C,\mathcal{O}_C(E))$ a global section with polar divisor $E$. 
		Then we have 
		\[ \tilde{\mu}_0(f\otimes s) \in H^0\textrm{\Large (}C, \omega_C\textrm{\large (}-E + \sum_{i=1}^m (\alpha_i -1)y_i + \sum_{j=1}^{n}\beta_jx_j\textrm{\large )}/\omega_C(-E)\textrm{\Large )} \]
		if and only if 
		\[ s\in H^0\textrm{\Large (}C, \omega_C\textrm{\large (}-2E + \sum_{i=1}^m (\alpha_i -1)y_i + \sum_{j=1}^{n}\beta_jx_j\textrm{\large )}\textrm{\Large )}\]
		
		Since $H^0(C, \mathcal{O}_C(E))$ is spanned by $1$ and $f$, it follows that our loci are transverse if and only if 
		\[ h^0\textrm{\Large (}C, \omega_C\textrm{\large (}-2E + \sum_{i=1}^m (\alpha_i -1)y_i + \sum_{j=1}^{n}\beta_jx_j\textrm{\large )}\textrm{\Large )} = 0\]
		From the Riemann-Roch theorem, this is equivalent to 
		\[ h^0(C, E+ \sum_{i=1}^m y_i + \sum_{k=1}^{g+1-d-m}z_k) = 2\]
		This follows immediately from Corollary 5 in \cite{Bud}. 
		Hence $\mathrm{Im}(\varphi_{\alpha,\beta})$ and $C^1_d$ are transverse and the number of pairs $(L,y_1,\ldots, y_m)$ with the desired properties can be expressed as $\varphi_{\alpha,\beta}^{*}(c^1_d)$, where $c^1_d$ is the class of the cycle $C^1_d$. 
		
		We know from \cite{ACGC1}, Chapter VII that 
		\[ c^1_d = \frac{\theta^{g-d+1}}{(g-d+1)!}- \frac{x\theta^{g-d}}{(g-d)!}\]
		where $x$ is the class of the divisor $X_p = \left\{ D \in C_d \ | \ D- p \geq 0 \right\}$ and $\theta$ is the pullback from the Jacobian of the theta divisor. 
		Finally we use the formulas in \cite{ACGC1}, Chapter VIII to deduce 
		\[ \varphi_{\alpha,\beta}^{*}(c^1_d) = \frac{g!}{d!(g+1-d-m)!}(2d+m-g-1-\sum_{i=1}^m\frac{1}{\alpha_i})\prod_{i=1}^m \alpha_i^2\]
		
		In the case $m=0$ the proposition still holds, using the convention that $\prod_{i=1}^{m} \alpha_i = 1$ and $\sum_{i=1}^m \frac{1}{\alpha_i} = 0$. 
	\end{proof}
	
	Next, we consider degree $d$ holomorphic maps $f\colon \mathbb{P}^1 \rightarrow \mathbb{P}^1$ with given ramification profiles $b_1, b_2$ and $b_3$ over the points $0,1,\infty$ and unramified elsewhere and ask what is their number up to isomorphism.  We denote this number by $N$ and count it for different choices of $b_1, b_2$ and $b_3$. 
	
	\begin{proposition} \label{table} We have the following table: 
		
		\begin{table}[h] \def\arraystretch{1.5}
			\centering 
			
			\begin{tabular}{|c|c|} 
				\hline
				Ramification profiles & Number $N$ of maps  \\ \hline 
				$b_1 = (m,n)$ with $m\neq n$ and $b_2=b_3 =(2,\ldots, 2)$ &  0 \\ \hline 
				$b_1 = (k,k), b_2 = b_3 = (2,\ldots, 2)$ & 1  \\ \hline  
				$b_1 = (2k), b_2=(4,2,\ldots, 2), b_3=(2,\ldots,2) $& $k-1$ \\ \hline
				$b_1 = (2k), b_2=(2,2,\ldots, 2), b_3=(2,\ldots,2,1,1) $& 1 \\ \hline
				$b_1 = (2k-1,1), b_2 = (4,2,\ldots, 2), b_3 = (2,\ldots, 2,1,1)$ & $k-1$ \\ \hline
				$b_1 = (2k-1,1), b_2=(3,2,2,\ldots, 2,1), b_3 =(2,\ldots,2)$ & 1 \\ \hline 
				$b_1 = (2k-1,1), b_2 = (4,2,\ldots, 2,1,1), b_3 = (2,\ldots, 2)$ & $k-2$ \\ \hline
				$b_1 = (m,n,1)$ with $|m-n|\neq 1$, $b_2 =(4,2,2,\ldots,2), b_3 = (2,\ldots, 2)$ & 0 \\ \hline
				$b_1 = (k,k-1,1), b_2 =(4,2,2,\ldots,2), b_3 = (2,\ldots, 2)$ & 1 \\ \hline	
				$b_1 = (2k+1), b_2 =b_3 =(2,\ldots,2,1)$ & 1 \\ \hline		
			\end{tabular}
		\end{table}
	\end{proposition}
	\begin{proof}
		Using Corollary 4.10 in \cite{Mir} this question translates into a purely combinatorial one: the number of such maps is equal to the number of conjugacy classes of 3-tuples $(\sigma_1,\sigma_2,\sigma_3)$ of permutations in $S_d$ having cycle types $b_1, b_2$ and $b_3$ respectively and satisfying $\sigma_1  \circ\sigma_2 =\sigma_3$. We give a proof for the third case in the table and claim that all other cases follow similarly. 
		
		Let us assume $b_1 = (2k), b_2 = (4,2,\ldots,2)$ and $b_3 = (2,\ldots, 2)$. We want to compute the number of solutions of $\sigma_1 \circ \sigma_2 = \sigma_3 $ up to conjugacy, where $\sigma_1, \sigma_2$ and $\sigma_3$ have cycle types $b_1, b_2$ and $b_3$. Since we are interested in solutions up to conjugacy, we can assume $\sigma_1 = (1,2,\ldots,2k)$. Moreover, we denote the 4-cycle in $\sigma_2$ by $(a,a+s,a+t,a+v)$. 
		
		The relation $\sigma_1\circ\sigma_2 =  \sigma_3$ implies that $\sigma_3$ must contain the transposition $(a,a+s+1)$, which in turn implies $\sigma_2(a+s+1) = a-1$. If $a-1$ is not in the 4-cycle, we further get $\sigma_3(a-1) = a+s+2$ and so on, until we ge get a condition for an element of the 4-cycle. 
		
		We split the problem into different cases depending on the order of the numbers $0 <s,t,v <2k$. 
		
		\textbf{Case I:} The order is $v<s<t$. Then we apply repeatedly the reasoning outlined previously and deduce $\sigma_2(a+s+t-s) = a+s-t$. It then follows that $t-s = 2k-v$ which implies the contradiction $s<v$. 
		
		\textbf{Case II:} The order is $s<t<v$. In this case, we obtain as in the first case that $t-s = 2k-v$. Moreover, the same argument implies $a+(a+v) = (a+s) + (a+t)$. As a consequence, $t=k$ and $v=k+s$. By reasoning repeatedly as explained, for every $s$ we get the unique solution 
		\[ \sigma_2 = (a, a+s, a+k, a+k+s) \prod_{i=1}^{k-s-1}(a-i,a+i+s) \prod_{i = k-s+1}^{k-1}(a-i,a+k+s+i) \]  
		\[ \sigma_3 = \prod_{i=0}^{k-s-1}(a-i,a+i+s+1) \prod_{i = k-s}^{k-1} (a-i, a+k+s+i+1) \]
		This solution is up to conjugation the one where $a=1$ and hence we get $k-1$ possible solutions depending on the value of $s$ between 1 and $k-1$. 
		
		Using the same method we conclude all other assumptions on the order yield no solution.
	\end{proof}
	
	\subsection{A Hurwitz space over the elliptic curves I} For an integer $k\geq 1$, we are interested in degree $2k$ maps from an elliptic curve to $\mathbb{P}^1$, having ramification profiles $b_1 =(2k)$, $b_2 =b_3 = (2,2,\ldots,2)$ and $b_4 = (2,1,1,\ldots,1)$ over four branch points $q_1,q_2,q_3$ and $q_4$. Let $B \coloneqq \left\{b_1,b_2,b_3,b_4\right\}$ and consider the Hurwitz scheme $\overline{H}_{2k,B}$. Following the method in \cite{AltCat} Section 5, we prove:
	\begin{proposition} \label{ellipticsix}
		For any $k\geq 1$, the map $\pi_k\colon\overline{H}_{2k,B}\rightarrow\overline{\mathcal{M}}_{1,1}$ remembering the point $p$ of total ramification over $q_1$ and stabilizing the source curve has degree 6. 
	\end{proposition} 
	\begin{proof}
		The proposition is clear for $k=1$ so we can assume $k\geq 2$. We consider $[E_\infty,p]$ the singular curve of $\overline{\mathcal{M}}_{1,1}$ and we compute the length of the cycle $\pi_k^{*}([E_\infty,p])$ which we know is equal to the degree of the map. 
		
		Let us denote by $\pi\colon X\rightarrow \Gamma$ an admissible cover mapped by $\pi_k$ to $[E_\infty,p]$. We denote by $R$ the rational component of $X$ mapping to $E_\infty$ and by $R_1$ the component collapsing to the node of the curve $E_\infty$. Furthermore, we will denote by $u$ and $v$ the two nodes where $R$ and $R_1$ are glued together. It follows from our notation that $R$ contains the totally ramified point $p$. As for the curve $\Gamma$, we denote by $\mathbb{P}_1$ the target of $R$, by $\mathbb{P}_2$ the target of $R_1$ and by $q$ the node. Finally, we denote by $f$ and $f_1$ the restriction of $\pi$ to $R$ and $R_1$ respectively. 
		
		In order to compute the length of the cycle, we distinguish three different cases for the admissible cover $\pi\colon X\rightarrow \Gamma$ in $\pi_k^{-1}([E_{\infty},p])$ depending on the position of the branch points on the components of $\Gamma$. For each such admissible cover we compute its multiplicity in $\pi_k^{*}([E_\infty,p])$.  
		
		\textbf{Case I:} The points $q_1$ and $q_4$ are on $\mathbb{P}_1$. In this case, we apply the Riemann-Hurwitz Theorem to $f$ and deduce that there are exactly two points above the node $q$. This implies that the degree of $f_1$ is $2k$. Hence the only two components of $X$ are $R$ and $R_1$. 
		
		Let $i$ and $j$ be the ramification orders of $f$ at the points $u$ and $v$. The first two rows in the table of Proposition \ref{table} imply that $i= j =k$ and there is a unique choice of the map $f_1$ up to the action of $PGL(2)$ on $R_1$ and  $\mathbb{P}_2$. We are now ready to describe the maps $f$ and $f_1$ up to the $PGL(2)$-action. 
		
		For $f$ we assume that $u=0,v=1$ and $p =\infty$. Then the map can be given as 
		\[f(t) = t^k(t-1)^k\]
		and the simple ramification point is $\frac{1}{2}$. The only non-trivial automorphism of $f$ is $\tau(t) = 1-t$, permuting the nodes $u$ and $v$. 
		
		For $f_1$ we first provide an implicit description of it. Consider $u = 0, v = \infty$ and let $\tau_1$ and $\tau_2$ be the automorphisms $\tau_1(t) = \frac{\xi}{t}$ and $\tau_2(t) =\xi^2t$ where $\xi$ is a primitive root of order $2k$. Let $G$ be the group of automorphisms generated by $\tau_1$ and $\tau_2$. Clearly it has order $2k$ and the map $f_1$ can be taken to be the quotient
		\[ f_1\colon \mathbb{P}^1\rightarrow \mathbb{P}^1/G \]
		The double ramification points are $\eta,\eta^5,\ldots, \eta^{4k-3}$ over a branch point and $\eta^3, \eta^7,\ldots, \eta^{4k-1}$ over another, where $\eta$ is a primitive root of order $4k$. Using this, we get an explicit description of such a map $f_1$ to be 
		\[ f_1(t) = \frac{(t^k-\eta^k)^2}{t^k}\] 
		It is simply checked that the morphisms of $f_1$ are the elements of $G$. 
		
		We hence found a unique point in $\overline{H}_{2k,B}$ over $[E_\infty,p]$ having an automorphism group of order $2k$. We want now to compute its multiplicity in $\pi_k^{*}([E_\infty,p])$.   
		
		For this, we look at the complete local ring of $[\pi\colon X\rightarrow\Gamma]$ in $\overline{H}_{2k,B}$. We know it is the ring of invariants of 
		\[ \mathbb{C}[[t_1,t_{1,1},t_{1,2}]]_{/(t_1 = t_{1,1}^k =t_{1,2}^k)}\]
		with respect to the group $\mathrm{Aut}_\pi(X)$ of automorphisms $\alpha$ of $X$ satisfying $\pi\circ \alpha = \pi$. 
		But this is equal to the ring of invariants of  
		\[ \mathbb{C}[[t_1,t_{1,1},t_{1,2}]]_{/(t_1 = t_{1,1}^k =t_{1,2}^k)}\]
		with respect to the action of the subgroup $\mathrm{Aut}_\pi^R(X) \leq \mathrm{Aut}_\pi(X)$ of automorphisms fixing the component $R$. This happens because after a suitable change of coordinates, we can assume that $\mathrm{Aut}_{\pi}(C)$ acts linearly on the parameter space $\Delta$. 
		
		In order to clarify this claim, let $\cup_{j=1}^k \Delta = \mathrm{Spec}(\mathbb{C}[[t_1,t_{1,1},t_{1,2}]]_{/(t_1 = t_{1,1}^k =t_{1,2}^k)})$ consisting of $k$ disks glued together at their respective origins and consider the universal deformation 
		\[
		\begin{tikzcd}
			\mathcal{C}  \arrow{r}{\pi} \arrow[swap]{dr}{} & \mathcal{P} \arrow{d}{} \\
			& \cup_{j=1}^k \Delta
		\end{tikzcd}
		\]
		where locally near the node $u$, the space $\mathcal{C}$ is given by $x_1\cdot y_1 = t_{1,1}$ and locally near the node $v$, the space $\mathcal{C}$ is given by $x_2\cdot y_2 = t_{1,2}$. Because the action of $\mathrm{Aut}_{\pi}(C)$ is linear on $\cup_{j=1}^k \Delta$ and moreover extends analytically to an action on $\mathcal{C}$, we see 
		that it is enough to understand how it acts at the level of coordinates of the central fibre $X$ of $\mathcal{C}$. In the standard coordinates $x_1= x_2 = y_1 =t$ at 0 and 1, and $y_2 =\frac{1}{t}$ at $\infty$ we have the following action of the automorphism group: 
		Consider the automorphism fixing $R$ and acting as multiplication by $\xi^{2j}$ on $R_1$. This automorphism fixes the coordinates $x_1$ and $x_2$ on $R$, multiplies $y_1$ by $\xi^{2j}$ and multiplies $y_2$ by $\xi^{-2j}$. It follows that this automorphism sends $t_{1,1}= x_1y_1$ to $\xi^{2j}t_{1,1}$ and $t_{1,2}$ to  $\xi^{-2j}t_{1,2}$. 
		In particular, the ring of invariants with respect to the subgroup  $\mathrm{Aut}_\pi^R(X)$ is 
		\[\mathbb{C}[[t_1,t_{1,1}t_{1,2}]]_{/(t_1 = t_{1,1}^k = t_{1,2}^k) } \cong \mathbb{C}[[t_1,t_{1,1}t_{1,2}]]_{/(t_1^2 = (t_{1,1}t_{1,2})^k)}\] 
		We apply the same reasoning to deduce that the automorphism acting as $\tau(t) = 1-t$ on $R$ and as $\tau_1(t) = \frac{\xi}{t}$ on $R_1$ sends $t_{1,1}$ to $-\xi t_{1,2}$ and $t_{1,2}$ to $-\xi^{-1}t_{1,1}$. It simply follows that the rings of invariants with respect to $\mathrm{Aut}_\pi^R(X)$ and $\mathrm{Aut}_\pi(X)$ are the same.


		Hence the local picture at such a point $[\pi\colon X\rightarrow \Gamma]$ is the following: 
		\[
		\begin{tikzcd}
			\cup_{j=1}^k \Delta \arrow{r}{k:1}  & (\cup_{j=1}^k \Delta)/\mathrm{Aut}_\pi(X) \arrow{d}{}  \subseteq \overline{H}_{2k,B} \\
			& \overline{\mathcal{M}}_{1,1}
		\end{tikzcd}
		\]
		
		The contribution of the point $[\pi\colon X\rightarrow \Gamma]$ to the length of $\pi_k^{*}([E_\infty,p])$ is equal to the local degree of $\pi_k$ at the point. We see from the diagram that this is equal to $\frac{1}{k}\cdot \mathrm{deg}(\cup_{j=1}^k \Delta \rightarrow \overline{\mathcal{M}}_{1,1})$. For each disk, the multiplicity of the point over $[\pi\colon X\rightarrow \Gamma]$ can be computed using curves in the universal deformation as exemplified in \cite{AltCat}, \cite{KodMg} and \cite{Diaz}.	The multiplicity of each of the $k$ points in $\amalg_{j=1}^k \Delta$ over $[E_\infty, p]$ is 2. It follows that the contribution of $[\pi\colon X\rightarrow \Gamma]$ to $\pi_k^{*}([E_\infty,p])$ is 2. 
		
		\textbf{Case II:} \ The points $q_1$ and $q_2$ are on $\mathbb{P}_1$. In this case, the Riemann-Hurwitz theorem applied to $f$ implies there are exactly $k+1$ points over $q$. It follows that $f_1$ is a degree 2 map and the ramified point above $q_4$ is on $R_1$. Furthermore, there are $k-1$ rational components $R_2, R_3,\ldots, R_k$ glued to $R$ at the $k-1$ points in $f^{-1}(q)\setminus \left\{u,v \right\}$. For every $j = \overline{2,k}$ the degree of the map $f_j = \pi_{|R_j}$ is equal to 2. 
		
		Observe that each $f_j$ has a unique non-trivial automorphism and moreover for $f_1$, this automorphism permutes $u$ and $v$.  
		
		We turn our attention to the map $f\colon R \rightarrow \mathbb{P}_1$. We know it has ramification types $(2k), (2,2,\ldots, 2)$ and $ (2,2,\ldots, 2,1,1)$ over three branch points and is otherwise unramified. 
		
		Up to conjugacy, there is a unique solution $\sigma_1 = \sigma_2 \circ \sigma_3$ in $S_{2k}$ for $\sigma_1, \sigma_2, \sigma_3$ permutations of cycle type $(2k), (2,2,\ldots, 2), (2,2,\ldots, 2,1,1)$, hence $f$ is unique up to the $\mathrm{PGL}(2)$-action on the curves $R$ and $\mathbb{P}_1$. 
		
		We will show that $f$ has a non-trivial automorphism. Such an automorphism is unique since it fixes $p$ and permutes $u$ and $v$. To see it exists, observe that $f$ can be written as a composition 
		\[ R \xrightarrow{2:1} \mathbb{P}^1 \xrightarrow{k:1} \mathbb{P}_1\]
		This can be deduced from unicity, reasoning on the parity of $k$. The $2:1$ map induces an involution on $R$, which is our desired automorphism. 
		
		Consider $\tau_1$ an automorphism of $X$ that is 
		non-trivial on the components $R$ and $R_1$ and permutes the components $R_2, \ldots, R_k$ accordingly. For $j = \overline{2,k}$ consider $\tau_j$ the automorphism of $X$ that is non-trivial on $R_j$ and trivial on all other components. Then the automorphism group $\mathrm{Aut}_\pi(X)$ is the group of cardinality $2^k$ generated by $\tau_1,\ldots, \tau_k$. 
		
		We know that the complete local ring of $[\pi\colon X \rightarrow \Gamma]$ in $\overline{H}_{2k,B}$ is the ring of invariants of 
		\[ \mathbb{C}[[t_1,t_{1,1}, t_{1,2},t_{1,3}\ldots,t_{1,k+1} ]]_{/(t_1 = t_{1,1} = t_{1,2} = t_{1,3}^2 = \cdots = t_{1,k+1}^2)} \] 
		with respect to $\mathrm{Aut}_\pi(X)$. The same reasoning as in the first case implies the complete local ring is isomorphic to
		\[ \mathbb{C}[[t_1]]\cong \mathbb{C}[[t_1,t_{1,1}\cdot t_{1,2},t^2_{1,3}\ldots,t^2_{1,k+1} ]]_{/(t_1 = t_{1,1} = t_{1,2} = t_{1,3}^2 = \cdots = t_{1,k+1}^2)} \] 
		
		Moreover for a local disk $\Delta_{t_1}$ around $[\pi\colon X\rightarrow \Gamma]$ we have a universal family 
		
		\[
		\begin{tikzcd}
			\mathcal{C}  \arrow{r}{\pi} \arrow[swap]{dr}{} & \mathcal{P} \arrow{d}{} \\
			& \Delta_{t_1}
		\end{tikzcd}
		\]
		with central fibre $\pi\colon X\rightarrow \Gamma$ and with local equations $x_iy_i = t_1$ at the nodes of $X$. We see that by collapsing the components $R_1,\ldots, R_k$ we get a family of genus 1 curves with central fiber $E_\infty$ and with local equation at the node of the form $xy =t_1^2$. Hence this point contributes with multiplicity 2 to the length of $\pi_k^{*}([E_\infty,p])$. 
		
		The third case when the points $q_1,q_2$ are in $\mathbb{P}_1$ is treated identically and we thus get another contribution of 2. To conclude, the length of $\pi_k^{*}([E_\infty,p]) = \mathrm{deg}(\pi_k)$ is $2+2+2  =6$. 
		
	\end{proof} 
	\subsection{A Hurwitz space over the elliptic curves II} Let again $k \geq 2$ and consider the partitions of $2k$: $b_1 =(2k-1,1)$, $b_2 =(4,2,\ldots,2), b_3 =(2,\ldots, 2)$ and $b_4=(2,1,1,\ldots,1)$. We denote by $B$ the set $\left\{b_1,b_2,b_3,b_4\right\}$ and we consider $\overline{H}_{2k,B}$ the Hurwitz scheme of admissible covers of degree $2k$ having ramification profiles $b_1,b_2,b_3,b_4$ over four points $q_1, q_2,q_3$ and $q_4$. We want to compute the degree of $\pi_k\colon \overline{H}_{2k,B} \rightarrow \overline{\mathcal{M}}_{1,1}$ remembering only the point $p$ of ramification order $2k-1$ and stabilizing the source curve. 
	\begin{proposition} \label{doiHurw} The degree of $\pi_k\colon\overline{H}_{2k,B} \rightarrow \overline{\mathcal{M}}_{1,1}$ is $6k-3$. 
	\end{proposition}
	\begin{proof} We consider the singular point $[E_\infty,p]$ of $\overline{\mathcal{M}}_{1,1}$ and we compute the length of $\pi_k^{*}([E_\infty,p])$, which we know is equal to the degree of the map. Our approach is again to consider all points in the preimage $\pi_k^{-1}([E_\infty,p])$ and compute their multiplicity.  
		
		Let $\pi\colon X\rightarrow \Gamma$ an admissible cover mapped by $\pi_k$ to $[E_\infty,p]$. In what follows, we preserve the notations in Proposition \ref{ellipticsix} for the components and nodes of the source and target. 
		
		Depending on the position of the branch points we distinguish again three different cases. 
		
		\textbf{Case I:} The points $q_1$ and $q_2$ are on $\mathbb{P}_1$. In this case, applying the Riemann-Roch theorem to $f\colon R \rightarrow \mathbb{P}_1$ we deduce there are $k+1$ points in the fibre over $q$. This is possible if and only if the ramification profile over $q$ is $(2,2,\ldots,2,1,1)$. Hence $f_1$ has degree 2 and the ramified point over $q_4$ is on $R_1$. Furthermore $X$ has another $k-1$ rational components $R_2,R_3,\ldots, R_k$ that are glued to $R$ at the $k-1$ points in $f^{-1}(q)\setminus \left\{u,v \right\}$. For every $j = \overline{2,k}$ the degree of the map $f_j = \pi_{|R_j}$ is equal to 2. 
		
		Such an admissible cover admits no non-trivial automorphism except the ones on the components $R_2,\ldots, R_k$. We see from Proposition \ref{table} that there are $k-1$ such admissible covers. Moreover, the complete local ring at such a point in $\overline{H}_{2k,B}$ is the ring of invariants of 
		\[ \mathbb{C}[[t_1,t_{1,1}, t_{1,2}, t_{1,3},\ldots, t_{1,k+1}]]_{/(t_1=t_{1,1} = t_{1,2} = t_{1,3}^2 =\cdots =t_{1,k+1}^2)}\]
		with respect to $\mathrm{Aut}_\pi(X)$. Using a similar argument as the one in Proposition \ref{ellipticsix}, we deduce that the ring of invariants is isomorphic to 
		\[ \mathbb{C}[[t_1]] \cong \mathbb{C}[[t_1,t_{1,1}, t_{1,2}, t^2_{1,3},\ldots, t^2_{1,k+1}]]_{/(t_1=t_{1,1} = t_{1,2} = t_{1,3}^2 =\cdots =t_{1,k+1}^2)} \]
		and we also deduce that each such admissible cover is counted with multiplicity $2$. Hence the admissible covers in this case contribute with $2k-2$ to the count. 
		
		\textbf{Case II:} The points $q_1$ and $q_3$ are on $\mathbb{P}_1$. Applying the Riemann-Hurwitz theorem to $f\colon R \rightarrow \mathbb{P}_1$ we deduce that the number of points in the fibre over $q$ is $k$. We have two possibilities: either the point of ramification order 4 is on the rational component collapsing to the node of $E_\infty$, or the point of ramification order 4 is on a rational component collapsing to a smooth point of $E_\infty$ when we stabilize the curve.
		
		For the first possibility, there is a unique choice for the admissible cover, as implied by the first and sixth rows in Proposition \ref{table}. The map $f\colon R \rightarrow \mathbb{P}_1$ has ramification profiles $(2k-1,1)$, $(2,\ldots, 2)$ and $(3,2,\ldots,2,1)$ over $q_1, q_3$ and $q$. For the map $f_1\colon R_1\rightarrow \mathbb{P}_2$ we know that the ramification profiles are $(4), (2,1,1)$ and $(3,1)$ over $q_2,q_4$ and $q$. At each of the other $k-2$ points over $q$, the curve $X$ contains a rational curve glued to $R$ that maps $2:1$ to $\mathbb{P}_2$ with branch points $q$ and $q_2$. 
		
		The complete local ring at this point is 
		\[ \mathbb{C}[[t_{1,1}]]\cong  \mathbb{C}[[t_1,t_{1,1}, t_{1,2}, t^2_{1,3},\ldots, t^2_{1,k}]]_{/(t_1=t^3_{1,1} = t_{1,2} = t_{1,3}^2 =\cdots =t_{1,k}^2)} \] 
		Hence at this point $\overline{H}_{2k,B}$ is smooth and, as the ramification orders at $u$ and $v$ are $3$ and $1$, the multiplicity of the point is $1+3 = 4$. 
		
		For the second possibility, Proposition \ref{table} implies there are $k-2$ choices of the map $f\colon R \rightarrow \mathbb{P}_1$. The map $f_1 \colon R_1 \rightarrow \mathbb{P}_2$ is $2:1$ with branch points $q_2$ and $q_4$. At the point of order 4 over $q$ we have a rational component glued to $R$ mapping $4:1$ to $\mathbb{P}_2$ with two points of total ramification over $q$ and $q_2$. At all the other points in the fibre over $q$ there is a rational component glued to $R$ that maps $2:1$ to $\mathbb{P}_2$ with branch points $q$ and $q_2$. 
		
		Again, we see that these points are smooth in $\overline{H}_{2k,B}$ and they are all counted with multiplicity $2$. Hence the contribution of Case II to the count is $4+2(k-2) = 2k$. 
		
		\textbf{Case III:} The points $q_1$ and $q_4$ are on $\mathbb{P}_1$. In this case the map $f\colon R \rightarrow \mathbb{P}_1$ can have degree $2k-1$ or $2k$. If we assume the degree in $2k$, the Riemann-Hurwitz theorem implies there is a unique point in the fibre over $q$, which is false, as both $u$ and $v$ are there. 
		
		It follows that the degree of $f$ is $2k-1$ and, by the Riemann-Hurwitz theorem, that $u$ and $v$ are the only points in the fibre over $q$. In this case, Proposition \ref{table} implies that $f_1\colon R_1\rightarrow \mathbb{P}_2$ is the unique map of degree $2k$ and ramification profile $(4,2,\ldots, 2), (2,\ldots,2), (k,k-1,1)$ over $q_2, q_3$ and $q$. In this case, the same argument with the complete local ring implies this cover should be counted with multiplicity $k-1+k = 2k-1$. 
		
		Adding up the three cases, it follows that the degree of $\pi_k$ is $2k-2+2k+2k-1 = 6k-3$. 
		
	\end{proof}
	\subsection{Combinatorial identities} Another ingredient we will require is the computation of some combinatorial sums. We state without proof the following identities 
	\begin{proposition} \label{combiden} We consider the sums $S_k = \sum_{s=0}^{i-1} s^k \binom{2i}{s}$ and compute the first terms to be 
		\[ S_0 = 2\cdot 2^{2i-2}- \binom{2i-1}{i}, \ S_1 = 2i\cdot 2^{2i-2} - 2i\cdot \binom{2i-1}{i} \ \textrm{and \ } S_2 = i(2i+1)\cdot 2^{2i-2} - 3i^2 \cdot \binom{2i-1}{i}\]
		Similarly, let $T_k = \sum_{s=0}^{i-1}s^k \binom{2i-1}{s}$. We compute the first terms to be
		\[T_0 = 2^{2i-2}, \ \ T_1 = \frac{2i-1}{2}\cdot2^{2i-2} - \frac{i}{2}\cdot \binom{2i-1}{i}, \ \ T_2 = \frac{(2i-1)i}{2}\cdot2^{2i-2}- \frac{(2i-1)i}{2}\cdot\binom{2i-1}{i}\]
		\[T_3 = (i^3-\frac{3}{4}i+\frac{1}{4})\cdot2^{2i-2}-(\frac{3}{2}i^3-i^2)\binom{2i-1}{i} \] 
		\[ \textrm{and} \ \ T_4 = (i^4+i^3-\frac{9}{4}i^2+\frac{3}{4}i)\cdot2^{2i-2}-(2i^4-i^3-i^2+\frac{1}{2}i)\binom{2i-1}{i}\]
	\end{proposition}
	
	\section{Test curves} \label{sec2}
	We consider some classical examples of test curves on $\overline{\mathcal{M}}_g$ and take their pullbacks to $\overline{\mathcal{R}}_g$ in order to obtain new ones. Using the projection formula and the description of $\overline{\mathcal{R}}_g$ over $\overline{\mathcal{M}}_g$ in \cite{Casa} and \cite{FarLud}, we compute their intersection numbers with all divisorial classes. All test curves on $\overline{\mathcal{M}}_g$ we are considering can be found in \cite{ModHM} and \cite{Mul}, while most test curves on $\overline{\mathcal{R}}_g$ we consider appear in more detail in \cite{Carlos}.
	
	Before providing the test curves we remark that $\mathrm{Pic}_\mathbb{Q}(\mathcal{R}_g) \cong \mathbb{Q}$. This follows from Theorem A and Theorem B in \cite{Putman}. Consequently, the classes $\lambda, \delta_0', \delta_0'', \delta_0^{\mathrm{ram}}$ together with $\delta_i, \delta_{g-i}, \delta_{i:g-i}$ for $1\leq i\leq [g/2]$ form a basis of  $\mathrm{Pic}_\mathbb{Q}(\overline{\mathcal{R}}_g)$. In particular, describing the intersection of the test curves with these classes is sufficient to describe their intersection with any class.
	\subsection{Test curve $A$} Let $A$ be the test curve in $\overline{\mathcal{M}}_g$ consisting of a generic genus $g-1$ curve $C$ glued at a generic point $x$ to a pencil of elliptic curves along a base point. Taking the pullback of the curve $A$ to $\overline{\mathcal{R}}_g$ we obtain three test curves $A_{g-1}, A_1, A_{1:g-1}$, contained in the three divisorial components $\Delta_{g-1}, \Delta_1$ and $\Delta_{1:g-1}$ respectively. We have the following intersection numbers, where the omitted intersections are all 0:
	\[ A_{g-1} \cdot \lambda = 1,\  A_{g-1} \cdot \delta_0' = 12, \ A_{g-1} \cdot \delta_{g-1} = -1 \]
	\[ A_1 \cdot \lambda = 3, \ A_1 \cdot \delta_0'' = 12, \ A_1 \cdot \delta^{\mathrm{ram}}_0 = 12, \ A_1 \cdot \delta_1 = -3  \]
	\[ A_{1:g-1} \cdot \lambda = 3, \ A_{1:g-1} \cdot \delta_0' = 12, \ A_{1: g-1} \cdot \delta^{\mathrm{ram}}_0 = 12, \ A_{1: g-1} \cdot \delta_{1: g-1} = -3 \]
	
	\subsection{Test curve $B$} Consider a generic point $[C,x] \in \mathcal{M}_{g-1,1}$. By glueing the point $x$ to a point $y$ moving on the curve, we obtain a curve $B$ on $\overline{\mathcal{M}}_g$. Proceeding as before, the pullback provides 3 test curves $B'$, $B''$ and $B^{\mathrm{ram}}$ lying in the divisors $\Delta_0'$, $\Delta''_0$ and $\Delta_0^{\mathrm{ram}}$. We have the following intersection numbers, the ones omitted being 0: 
	\[ B' \cdot \delta_0' = (1-g)(2^{2g}-4), \ B' \cdot \delta_{g-1} = 2^{2g-2}-1, \ B' \cdot \delta_{1: g-1 } =  2^{2g-2}-1 \]
	\[ B'' \cdot \delta_1 = 1, \  B'' \cdot \delta_0'' = 2-2g  \] 
	\[ B^{\mathrm{ram}} \cdot \delta_0^{\mathrm{ram}} = 2^{2g-2}(1-g), \ B^{\mathrm{ram}} \cdot \delta_1 = 1, \ B^{\mathrm{ram}} \cdot \delta_{1: g-1} = 2^{2g-2}-1 \] 
	
	\subsection{Test curves $C_i$} Let $i$ an integer satisfying $2\leq i \leq g-1$ and let $[C] \in \mathcal{M}_i$ and $[D,y] \in \mathcal{M}_{g-i,1} $ be two generic curves. Let $\eta_C \in \mathrm{Pic}(C)[2]\setminus \left\{0\right\}$ and $\eta_D \in $$\mathrm{Pic}(D)[2]\setminus \left\{0\right\}$ and consider the test curves in $\overline{\mathcal{R}}_g$ 
	\[[C\cup_{x \sim y } D, (\eta_C, \mathcal{O}_D)]_{x\in C} \]  
	\[[C\cup_{x \sim y } D, (\mathcal{O}_C, \eta_D)]_{x\in C} \]  
	\[[C\cup_{x \sim y } D, (\eta_C, \eta_D)]_{x\in C} \] 
	by varying $x$ along $C$. We denote them $C^i_{i}, C^i_{g-i} $ and $C^i_{i: g-i}$ respectively. It is clear they are  contained in the divisors $\Delta_i, \Delta_{g-i}$ and $\Delta_{i:g-i}$ respectively. The intersection numbers are the following, where all omitted intersection numbers are 0: 
	\[ C^i_i \cdot \delta_i = 2-2i \]
	\[ C^i_{g-i} \cdot \delta_{g-i} = 2-2i \]
	\[ C^i_{i: g-i} \cdot \delta_{i\colon g-i} = 2-2i \] 
	\section{Intersection numbers}
	
	Throughout this section, the genus $g = 2i$ is even and the partition $\mu$ of length $g-1$ is taken to be $(4,2,2,\ldots, 2, -2,\ldots,-2)$. Our goal is to compute the intersection of $\overline{D}(\mu)$ with some of the test curves described in Section \ref{sec2} and conclude Theorem \ref{maintheorem}.  
	
	We consider the normalization $\nu\colon \overline{H}^\nu_{g,\mu} \rightarrow\overline{H}_{g,\mu}$ of the Hurwitz scheme $\overline{H}_{g,\mu}$. Our first task is to extend the map $c_\mu\circ\nu\colon \overline{H}^\nu_{g,\mu} \dashrightarrow \overline{\mathcal{R}}_g$ over points in the preimage $(\pi_\mu\circ\nu)^{-1}([C_{/x\sim y}])$ where $[C,x]$ is a generic point of $\mathcal{M}_{g-1,1}$. Our approach is to describe the admissible covers in $\overline{H}_{g,\mu}$ above such $[C_{/x\sim y}]$ in $\overline{\mathcal{M}}_g$ and based on the description explain how the map $c_\mu\circ \nu$ is extended. 
	
	\begin{proposition} \label{extension} Let $[C,x]$ be a generic point of $\mathcal{M}_{g-1,1}$ and $y$ a point on the curve $C$. Then the rational map $c_\mu\circ\nu\colon \overline{H}^\nu_{g,\mu} \dashrightarrow \overline{\mathcal{R}}_g$ can be extended over the fibre $(\pi_\mu\circ\nu)^{-1}([C_{/x\sim y}])$.
	\end{proposition}
	\begin{proof} Consider an element $[\pi\colon X\rightarrow \Gamma]$ in $\overline{H}_{g,\mu}$ above  $[C_{/x\sim y}]\in\overline{\mathcal{M}}_g$. The genericity of $[C,x]$ and Lemma 5.3 in \cite{Diaz} imply that $\Gamma$ has exactly two irreducible components. We denote $\Gamma = \mathbb{P}_1 \cup_q \mathbb{P}_2$ where $C$, seen as a component of $X$, is mapped to $\mathbb{P}_1$. We distinguish four different possibilities for the admissible cover $\pi\colon X\rightarrow \Gamma$ depending on the position of the points $q_1$ and $q_2$ on $\Gamma$. 
		
		\textbf{Case I:} The points $q_1$ and $q_2$ are on $\mathbb{P}_1$. In this case, it follows that $\mathrm{deg}(\pi_{|C}) = g$. Otherwise, the twist on $C$ will correspond to a partition of length strictly smaller than $g-1$, contradicting the genericity of the curve. Moreover, as the map 
		\[ \mathcal{H}^0_{g-1}(2,1,1,\ldots,1, -1,\ldots,-1) \rightarrow \mathcal{M}_{g-1}\]
		corresponding to the length $g-1$ partition $(2,1,1,\ldots,1, -1,\ldots,-1)$ is not dominant, it follows that 
		\[ \mathcal{O}_C(\frac{\pi^{*}(q_1)-\pi^{*}(q_2)}{2}) \not\cong \mathcal{O}_C \]
		Taking into account the genericity of the point $x$, we deduce that the curve $X$ is as in the following figure 
		\begin{figure}[H] \centering 
			\begin{tikzpicture}
				\draw [domain=0.3:2] plot ({\x},{sqrt{4*\x-\x*\x} - 2*\x});
				\draw [domain=-2:5] plot ({\x*\x/5-\x/2+1.3)},{\x/5-3.6});
				\draw [domain=0.6:2.3] plot ({\x}, {\x/15-2});
				\draw [domain=0.6:2.3] plot ({\x-0.5}, {\x/15-1}); 
				
				\node[left = 1mm of {(1.3,0)}] {$C$};
				\node[left = 1mm of {(3.3,-3)}] {$R$};
				
				\fill (0.95,-0.9)  circle[radius=1.5pt];
				\fill (1.37,-1.9)  circle[radius=1.5pt];
				\fill (1.72,-2.97) circle[radius=1.5pt];
				\node[right = 1mm of {(1.65,-2.8)}] {$x$};
				\fill (1.95,-3.8) circle[radius=1.5pt];
				\node[right = 1mm of {(1.8,-3.6)}] {$y$};
				\draw [->]  (2.2, -1.4) -- (1.5,-1.8) ;
				\draw [->]  (2.2, -1.4) -- (1.1,-1) ;
				\node[right = 1mm of {(-0.3,-0.75)}] {$\text{\small 1:1}$};
				\node[right = 1mm of {(0,-1.7)}] {$\text{\small 1:1}$};
				\node[right = 1mm of {(2.2,-1.4)}] {The other $g-2$ points in the same fiber as $x$};
				\node[right = 1mm of {(0.4,-3)}] {$2:1$};
			\end{tikzpicture} \caption{A curve $X$ corresponding to an admissible cover in $\overline{H}_{g,\mu}$ over $[C_{/x\sim y}]$ } \label{figure}
		\end{figure}
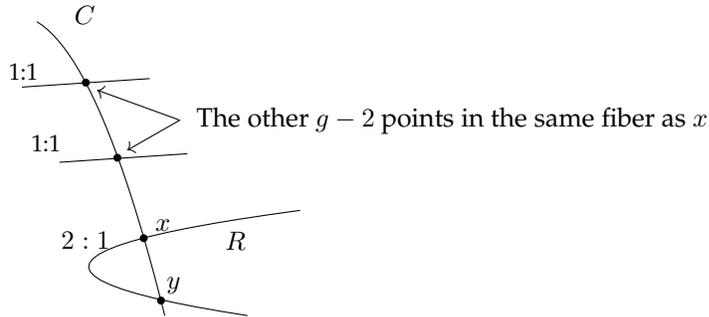 
		Here $y$ is one of the other $g-1$ points in the same fiber of $\pi$ as $x$ and there is a rational component mapping $1:1$ to $\mathbb{P}_2$ glued at the other $g-2$ points of this fibre. The points $x$ and $y$ are glueing $C$ to a rational component $R$ mapping $2:1$ to $\mathbb{P}_2$. 
		
		As there exists no non-trivial automorphism $\alpha\colon X\rightarrow X$ satisfying $\pi\circ\alpha = \pi$, the complete local ring of $[\pi\colon X \rightarrow \Gamma]$ in $\overline{H}_{g,\mu}$ has the form 
		\[ \mathbb{C}[[t_{1,1}, \ldots, t_{1,g}, t_1,\ldots, t_{3g-4}]]_{/(t_1 = t_{1,1} = \cdots = t_{1,g})} \cong \mathbb{C}[[t_1,\ldots, t_{3g-4}]]\]
		
		This implies that $\overline{H}_{g,\mu}$ is smooth at $[\pi\colon X\rightarrow \Gamma]$. We then consider an open $(3g-4)$-dimensional polydisk $\Delta^{3g-4}$ centered at $[\pi\colon X \rightarrow \Gamma]$ and assume the equation $t_1 = 0$ parametrizes the locus of admissible covers with singular source curve. We have the commutative diagram of maps 
		
		\[
		\begin{tikzcd}
			& \overline{\mathcal{R}}_g \arrow{d}{} \\
			\Delta^{3g-4} \setminus \left\{t_1 = 0\right\} \arrow{r}{\pi_\mu} \arrow[swap]{ur}{c_\mu} & \overline{\mathcal{M}}_g 
		\end{tikzcd}
		\] 
		
		Consider a small enough neighbourhood $U$ of $[C_{/x\sim y}]$ in $\overline{\mathcal{M}}_g$ such that its preimage in $\overline{\mathcal{R}}_g$ consists of disjoint open sets, each containing a unique element in the fibre over $[C_{/x\sim y}]$. We further shrink $\Delta^{3g-4}$ so that the image is contained in $U$. It follows that $\Delta^{3g-4} \setminus \left\{t_1 =0\right\}$ is mapped by $c_\mu$ to a unique open component of the preimage of $U$ in $\overline{\mathcal{R}}_g$ and the map $c_\mu$ can be analytically extended due to Hartogs' extension theorem.
		
		We want to determine the boundary divisor of $\overline{\mathcal{R}}_g$ to which $[\pi\colon X \rightarrow \Gamma]$ is mapped by this extension. For this we consider a 1-dimensional smoothing of the admissible cover 
		\[
		\begin{tikzcd}
			\mathcal{C}  \arrow{r}{\pi} \arrow[swap]{dr}{} & \mathcal{P} \arrow{d}{} \\
			& \Delta  
		\end{tikzcd}
		\]
		together with the sections $q_i \colon \Delta \rightarrow \mathcal{P}$ defining the branch points. By eventually shrinking $\Delta$ we can assume that the preimages $\pi^{-1}(q_1(\Delta))$ and $\pi^{-1}(q_2(\Delta))$ are set-theoretically unions of disjoint sections of $\mathcal{C}\rightarrow \Delta$. We consider the divisor $Z$ on $\mathcal{C}$ satisfying $2Z= \pi^{*}(q_{1*}(\Delta)) - \pi^{*}(q_{2*}(\Delta))$ and take the associated line bundle $\mathcal{O}_\mathcal{C}(Z)$. We then have that 
		\[ \mathcal{O}_\mathcal{C}(Z)_{|C_t} = \mathcal{O}_{C_t}(Z_t) \in \mathrm{Pic}(C_t)[2] \textrm{ \ for all \ } t \neq 0 \]
		
		As $\mathcal{C}$ is smooth, it follows that there exist coefficients $c, c_1, \ldots, c_{g-2}$ such that 
		\[ \mathcal{O}_\mathcal{C}(2Z) \cong \mathcal{O}_\mathcal{C}(cR + \sum_{j=1}^{g-2}c_jR_j) \]
		where $R, R_1,\ldots, R_{g-2}$ are the rational components of the central fibre $X$. We know that these rational components are disjoint with $Z$ and the intersection numbers of the rational components are $R \cdot R = -2$, $R\cdot R_j = 0 \ \forall j = \overline{1,g-2}$ and $R_j \cdot R_j = -1 \ \forall j = \overline{1,g-2}$. 
		As a consequence, it follows that $c = c_1 = \cdots =c_{g-2}$ = 0. Hence the central fiber of $\mathcal{O}_\mathcal{C}(Z)$ is a $2$-torsion point of $\mathrm{Pic}(X)$ and this line bundle is trivial on all the rational components. Hence, by collapsing the rational components we obtain an element of $\overline{\mathcal{R}}_g$ over the central fibre. Moreover, this element is in the divisor $\Delta_0'$ of $\overline{\mathcal{R}}_g$  because we have 
		\[ \mathcal{O}_C(\frac{\pi^{*}(q_1)-\pi^{*}(q_2)}{2}) \not\cong \mathcal{O}_C \]
		
		\textbf{Case II:} The points $q_1$ and $q_2$ are on $\mathbb{P}_2$. In order to treat this case, we need to introduce some notations. Denote by $x_1,\ldots, x_s$ the points on $C$ in the fibre $\pi^{-1}(q)$ that are different  from $x$ and $y$. For each rational component $\mathbb{P}$ of $X$ we are interested to which component of $\Gamma$ it is mapped and also, to which of the points $x\sim y, x_1,\ldots, x_s$ it collapses when we stabilize $X$. We introduce the following notations
		\[ c_{0,1} = \sum_{\substack{\mathbb{P} \textrm{ collapses to $x\sim y$} \\ \mathbb{P} \textrm{ is mapped to } \mathbb{P}_1}} \textrm{deg}(f_{|\mathbb{P}}\colon\mathbb{P} \rightarrow \mathbb{P}_1) \] 
		\[ c_{0,2} = \sum_{\substack{\mathbb{P} \textrm{ collapses to $x\sim y$} \\ \mathbb{P} \textrm{ is mapped to } \mathbb{P}_2}} \textrm{deg}(f_{|\mathbb{P}}\colon\mathbb{P} \rightarrow \mathbb{P}_2) \]
		
		We similarly define $c_{j,1}$ and $c_{j,2}$ for the points $x_j$ for $j = \overline{1,s}$. By double counting the ramification orders at the nodes we deduce the following numerical conditions for the admissible cover $\pi$: 
		\[ \mathrm{ord}_x(\pi_{|C}) + \mathrm{ord}_y(\pi_{|C})+c_{0,1} = c_{0,2}\]
		\[ \mathrm{ord}_{x_1}(\pi_{|C}) +c_{1,1} = c_{1,2} \]
		\[\ldots\]
		\[ \mathrm{ord}_{x_s}(\pi_{|C}) +c_{s,1} = c_{s,2} \] 
		
		We claim that $c_{j,1} \geq 1$ for all $j =\overline{1,s}$ and moreover that $c_{j,1} \geq 2$ if the point of ramification order $4$ is on a rational component collapsing to $x_j$. We also have a description for when equality holds. When $c_{j,1} = 1$ the only possibility is that $c_{j,2} = 2$. Moreover if $c_{j,1} =2$ and the point of ramification order $4$ collapses to $x_j$ it follows that $c_{j,2} = 4$. 
		
		To see this, assume that $c_{j,1} = 0$ for some $j$. This implies that there exists a unique rational component $\mathbb{P}$ that collapses to $x_j$. This component is furthermore mapped by $\pi$ to $\mathbb{P}_2$ and is totally ramified at the point of $\mathbb{P}$ glued to $x_j$. But the map $\pi_{|\mathbb{P}}$ has only even ramification orders over $q_1$ and $q_2$. By the Riemann-Hurwitz theorem, this is impossible. If $c_{j,1} =1 $ it follows again that there is a unique component $\mathbb{P}$ collapsing to $x_j$ and mapped to $\mathbb{P}_2$. This component $\mathbb{P}$ has ramification profiles $(2k-1,1), (2,\ldots, 2)$ and $(2,\ldots, 2)$ over $q, q_1$ and $q_2$. Reasoning as in Proposition \ref{table}, this is only possible if $k=1$ and hence $c_{j,2} =2$. The case when the point of ramification order $4$ is on a rational component collapsing to $x_j$ follows analogously. 
		
		We distinguish two different cases depending on whether this point of ramification order $4$ is collapsing to $x\sim y$ or to one of the points $x_1, \ldots, x_s$. 
		
		If this point collapses to $x \sim y$ we get by adding the equalities that 
		\[\mathrm{deg}(\pi_{|C}) + s + c_{0,1} \leq  2i \] 
		Furthermore the genericity of $[C,x]$ implies that 
		\[\mathrm{deg}(\pi_{|C}) + s \geq g-1 = 2i-1 \]
		
		Because $c_{0,1} \leq 1$ it follows that there is a unique rational component $R$ mapping to $\mathbb{P}_2$ that collapses to $x \sim y$. The Riemann-Hurwitz theorem applied to $\pi_{|R}$ implies that the fibre $\pi_{|R}^{-1}(q)$ contains at least $3$ distinct points. It follows that $c_{0,1} =1, \mathrm{deg}(\pi_{|C}) + s = 2i-1$ and the ramification profiles of $\pi_{|R}$ over the branch points $q, q_1$ and $q_2$ are $(m,n,1), (4,2,2,\ldots, 2)$ and $(2,\ldots, 2)$ where the entries add up to $2i-2s$.  
		
		Proposition \ref{table} implies that $m = i-s, n = i-s-1$ and in this case there is a unique choice for the map $\pi_{|R}$. We get two possible types of admissible covers $\pi\colon X \rightarrow \Gamma$ which we will now describe. 
		
		For the first type, $\pi_{|C}$ has degree $2i-s-1$ and the ramification orders at $x$ and $y$ are $i-s$ and $i-s-1$. Furthermore $\pi_{|R}$ is the unique map with ramification profiles $(i-s,i-s-1,1), (4,2,2, \ldots, 2)$ and $(2,2,\ldots, 2)$ over $q, q_1$ and $q_2$, where the points of ramification orders $i-s$ and $i-s-1$ are $x$ and respectively $y$. 
		
		For the second type, $\pi_{|C}$ has degree $2i-s-1$ and the ramification orders at $x$ and $y$ are $i-s-1$ and $i-s$. Furthermore $\pi_{|R}$ is the unique map with ramification profiles $(i-s,i-s-1,1), (4,2,2, \ldots, 2)$ and $(2,2,\ldots, 2)$ over $q, q_1$ and $q_2$, where the points of ramification orders $i-s-1$ and $i-s$ are $x$ and respectively $y$. 
		
		We see that for both types, there exists no non-trivial automorphism $\alpha\colon X \rightarrow X$ which satisfies $\pi \circ \alpha = \pi$. As a consequence the complete local ring of $\overline{H}_{g,\mu}$ at the point $[\pi\colon X\rightarrow \Gamma]$ is isomorphic to 
		\[\mathbb{C}[[t_{1,1},t_{1,2},t_2,\ldots, t_{3g-4}]]_{/(t^{i-s}_{1,1} = t_{1,2}^{i-s-1}) } \]  
		By considering $[\pi\colon X \rightarrow \Gamma]$ as a point in the normalization $\overline{H}_{g,\mu}^{\nu}$ of $\overline{H}_{g,\mu}$  we see that $[\pi\colon X \rightarrow \Gamma]$ is a smooth point of $\overline{H}_{g,\mu}^{\nu}$. The same argument using Hartogs' extension theorem as in Case I implies that the map $c_\mu\circ\nu$ can be extended over this point. As in Case I, we deduce that $[\pi\colon X\rightarrow\Gamma]$ is mapped to $\Delta''_0$ in $\overline{\mathcal{R}}_g$. 
		
		We now treat the case when the point of order $4$ is on a rational component collapsing to a smooth point of $C$. In this case we get by adding the equalities that  
		\[ \mathrm{deg}(\pi_{|C})+s +1 + c_{0,1} \leq 2i \]
		
		Using the inequality 
		\[ \mathrm{deg}(\pi_{|C}) + s\geq 2i-1\]
		implied by the genericity of $[C,x]$ we deduce that $\mathrm{deg}(\pi_{|C}) = 2i-s-1$ and $c_{0,1}= 0$. It follows that $\mathrm{ord}_x(\pi_{|C}) + \mathrm{ord}_y(\pi_{|C}) = 2i-2s-2$. Proposition \ref{table} implies that the ramification profiles of $\pi_{|R}$ over $q, q_1$ and $q_2$ are $(i-s-1,i-s-1), (2,\ldots,2)$ and $(2,\ldots, 2)$. 
		
		We see that the admissible cover $\pi\colon X \rightarrow \Gamma$ can be described as follows. The map $\pi_{|C}$ has degree $2i-1-s$ containing in its fiber over $q$ the points $x,y$ of ramification order $i-s-1$ and another point of ramification order 2 to which the point of ramification order $4$ over $q_1$ collapses. For the rational component $R$, the map $\pi_{|R}$ has ramification profiles $(i-s-1, i-s-1), (2,\ldots,2)$ and $(2,\ldots,2)$ over $q, q_1$ and $q_2$. 
		
		As a consequence, the complete local ring of $\overline{H}_{g,\mu}$ at $[\pi\colon X \rightarrow \Gamma]$ is the ring of invariants of 
		\[ \mathbb{C}[[t_{1,1}, t_{1,2}, t_{1,3}, t_2,\ldots, t_{3g-4}]]_{/(t_{1,1}^{i-s-1} = t_{1,2}^{i-s-1} = t_{1,3}^2)} \] 
		with respect to the group $\mathrm{Aut}_{\pi}(X)$ of automorphisms $\alpha\colon X\rightarrow X$ satisfying $\pi \circ \alpha = \pi$. 
		
		Denote by $x_1$ the point of ramification order $2$ in the same fibre as $x$ and $y$ and by $R_1$ the rational component of degree $4$ over $\mathbb{P}_2$ collapsing to $x_1$. Up to the $\mathrm{PGL}(2)$-action on both the source and the target we can describe $\pi_{|R_1}$ as 
		\[ \pi_{|R_1}(t) = t^2(t-1)^2\]
		In this description, observe that the point $x_1$ corresponds to the point $\frac{1}{2}$ and there is a non-trivial automorphism $\tau$ given as $\tau(t) = 1-t$. 
		
		By first considering the ring of invariants with respect to the automorphism acting as $\tau$ on $R_1$ and fixing $C$ and $R$ we get that the complete local ring of $\overline{H}_{g,\mu}$ at $[\pi\colon X \rightarrow \Gamma]$ is the ring of invariants of 
		\[ \mathbb{C}[[t_{1,1}, t_{1,2}, t_{1,3}^2, t_2,\ldots, t_{3g-4}]]_{/(t_{1,1}^{i-s-1} = t_{1,2}^{i-s-1} = t_{1,3}^2)} \cong \mathbb{C}[[t_{1,1}, t_{1,2}, t_2,\ldots, t_{3g-4}]]_{/(t_{1,1}^{i-s-1} = t_{1,2}^{i-s-1})}\]
		with respect to the automorphism subgroup $\mathrm{Aut}_\pi^{R_1}(X) \leq \mathrm{Aut}_\pi(X)$ of elements restricting to the identity on $R_1$. But $\mathrm{Aut}_\pi^{R_1}(X)$ is the group of automorphisms of $R$ fixing the points $x,y$ and all components of $X$ different from $R$. A description of this group appears in the proof of Proposition \ref{ellipticsix}. We consider the parameter space of 
		\[\mathbb{C}[[t_{1,1}, t_{1,2}, t_2,\ldots, t_{3g-4}]]_{/(t_{1,1}^{i-s-1} = t_{1,2}^{i-s-1})}\]
		consisting of $i-s-1$ polydisks of dimension $3g-4$ glued together along $t_{1,1} = 0$. Then $\mathrm{Aut}_\pi^{R_1}(X)$ acts on this space by identifying the $i-s-1$ components with one another. As a consequence, the space $\overline{H}_{g,\mu}$ is smooth at $[\pi\colon X \rightarrow \Gamma]$ and the method of Case I applies again. By choosing a 1-dimensional smoothing we deduce that $[\pi\colon X\rightarrow \Gamma]$ is mapped by $c_\mu$ to $\Delta_0''$ in $\overline{\mathcal{R}}_g$. 
		
		\textbf{Case III:} The point $q_1$ is on $\mathbb{P}_1$ and the point $q_2$ is on $\mathbb{P}_2$. The image $\mathrm{StMd}(X, \pi^{-1}(q_1), \pi^{-1}(q_2))$ by $a_\mu$ is an element of $\overline{\mathcal{H}}^0_g(\mu)/\mathcal{S}(\mu)$, hence there is a twist associated to this pointed curve. The genericity of $[C,x]$ implies that the partition on $C$ determined by the twist has length $g$. This imply that on $C$ we have a divisorial equivalence of the form 
		\[ 4y_1 + 2y_2 + \cdots + 2y_{i-1}-2x_1-\cdots -2x_{i-1}-x-y \equiv 0 \]
		in $\mathrm{Pic}(C)$ where $y_1, \ldots, y_{i-1},x_1,\ldots, x_{i-1}$ are points of $C$. In particular, this description implies $\mathrm{deg}(\pi_{|C}) = g, \mathrm{deg}(\pi_{|R}) = 2$ and at every point $x_j$ there is a rational component $R_j$ glued to $C$ that maps $2:1$ to $\mathbb{P}_2$ with branch points $q$ and $q_2$. 
		
		From the above description it follows that the complete local ring of $[\pi\colon X\rightarrow \Gamma]$ in $\overline{H}_{g,\mu}$ is the ring of invariants of 
		\[ \mathbb{C}[[t_{1,1},t_{1,2}, \ldots, t_{1,i+1}, t_2,\ldots, t_{3g-4}]]_{/(t_{1,1}^2 = \cdots = t_{1,i-1}^2 = t_{1,i}= t_{1,i+1})} \]
		with respect to the action of $\mathrm{Aut}_\pi(X)$. But $\mathrm{Aut}_\pi(X)$ is the group of cardinality $2^{i-1}$ generated by the automorphisms $\tau_1,\ldots,\tau_{i-1}$ where $\tau_j$ for $j = \overline{1,i-1}$ is the automorphism of $X$ that restricts to the identity on all components except $R_j$. In particular, the complete local ring is isomorphic to 
		\[ \mathbb{C}[[t_{1,1}^2,\ldots, t_{1,i-1}^2, t_{1,i}, t_{1,i+1}, t_2,\ldots, t_{3g-4}]]_{/(t_{1,1}^2 = \cdots = t_{1,i-1}^2 = t_{1,i}= t_{1,i+1})} \cong \mathbb{C}[[t_{1,i},t_2,\ldots, t_{3g-4}]] \]
		hence the point is smooth and the map $c_\mu$ can be extended to this point by Hartogs' extension theorem. 
		
		Next, we consider a 1-dimensional smoothing of $\pi\colon X \rightarrow \Gamma$ 
		\[
		\begin{tikzcd}
			\mathcal{C}  \arrow{r}{\pi} \arrow[swap]{dr}{} & \mathcal{P} \arrow{d}{} \\
			& \Delta  
		\end{tikzcd}
		\]
		together with the sections $q_1, q_2 \colon \Delta \rightarrow \mathcal{P}$. We consider as in Case I the divisor $Z$ on $\mathcal{C}$ satisfying $2Z = \pi^{*}(q_{1*}(\Delta)) - \pi^{*}(q_{2*}(\Delta))$. We see that 
		\[\mathrm{deg}(\mathcal{O}_\mathcal{C}(Z)_{|R}) = -1, \ \mathrm{deg}(\mathcal{O}_\mathcal{C}(Z)_{|C}) = i \mathrm{\ and \ } \mathrm{deg}(\mathcal{O}_\mathcal{C}(Z)_{|R_j}) = -1 \textrm{ \ for all \ } j = \overline{1,i-1}\]
		
		Because we have the self-intersection numbers $R\cdot R = -2$ and $R_j\cdot R_j = -1$ we conclude that by twisting with $-R- \sum_{j=1}^{i-1}R_j$ and collapsing the rational components $R_1,\ldots, R_{i-1}$ we obtain an element of $\Delta_0^{\mathrm{ram}}$ in $\overline{\mathcal{R}}_g$. 
		
		Indeed, when restricting $\mathcal{O}_\mathcal{C}(Z-R-\sum_{j=1}^{i-1}R_j)$ to $C$ we obtain the line bundle  
		\[L\cong \mathcal{O}_C(2y_1+y_2+\cdots+y_{i-1}-x_1-\cdots -x_{i-1}-x-y)\]
		and it is clear that $L^2 \cong \mathcal{O}_C(-x-y)$. Moreover the degree of $\mathcal{O}_\mathcal{C}(Z-R-\sum_{j=1}^{i-1}R_j)$ restricted to $R$ is 1, hence we obtain an element in $\Delta^{\mathrm{ram}}_0$ over the central fibre, as stated. 
		
		\textbf{Case IV:} The point $q_1$ is on $\mathbb{P}_2$ and the point $q_2$ is on $\mathbb{P}_1$. The image $\mathrm{StMd}(X,\pi^{-1}(q_1),\pi^{-1}(q_2))$ through $a_\mu$ is an element of $\overline{\mathcal{H}}^0_g(\mu)/\mathcal{S}(\mu)$. Using again the genericity of $[C,x]$ and the existence of a twist we deduce that on $C$ the twist determines a linear equivalence of one of the following forms
		\[ 4x_1+2x_2 + \cdots +2x_{i-2}+x+y - 2y_1-\cdots - 2y_i \equiv 0 \]
		\[ 2x_2+\cdots + 2x_{i-1} + (4-k)x + ky -2y_1-\cdots - 2y_i \equiv 0\]
		where all $x_j$ and $y_j$ are points of $C$ and $k \in \left\{1,2,3\right\}$. 
		
		We claim that the map $c_\mu\circ\nu\colon \overline{H}_{g,\mu}^\nu \rightarrow \overline{\mathcal{R}}_g$ extends over the preimages of such admissible covers $[\pi\colon X \rightarrow \Gamma]$ and the image of such a point is in $\Delta_0'$ only if the linear equivalence determined by the twist is of the second form and $k=2$. If the linear equivalence is of any of the other three forms, $c_\mu$ maps such admissible covers to $\Delta^{\mathrm{ram}}_0$. 
		
		The case when the divisorial equivalence is of the first form follows analogously to Case III. If the divisorial equivalence is of the second form with $k=2$, we see that the complete local ring of $\overline{H}_{g,\mu}$ at the point $[\pi\colon X\rightarrow \Gamma]$ is isomorphic to 
		\[\mathbb{C}[[t_{1,1}, t_{1,2}, t_2,\ldots, t_{3g-4}]]_{/(t_{1,1}^2 = t_{1,2}^2)}\]
		Hence the parameter space consists of two $(3g-4)$-dimensional polydisks glued together along the locus $t_{1,1}=0$ where the coordinates of the polydisks are considered to be $t_{1,1}, t_2,\ldots, t_{3g-4}$. In particular, both preimages of $[\pi\colon X \rightarrow \Gamma]$ are smooth in the normalization $\overline{H}_{g,\mu}^\nu$. 
		
		Over both polydisks we have a universal covering and hence the same approach as in Case I can be applied to deduce that $c_\mu\circ\nu$ maps the preimages of $[\pi\colon X \rightarrow \Gamma]$ to $\Delta_0'$. 
		
		We are left to treat the case when the divisorial equivalence is of the second form and $k=1$ or $k=3$. In both cases, the complete local ring of $\overline{H}_{g,\mu}$ at the point $[\pi\colon X\rightarrow \Gamma]$ is isomorphic to 
		\[\mathbb{C}[[t_{1,1}, t_{1,2}, t_2,\ldots, t_{3g-4}]]_{/(t_{1,1}^3 = t_{1,2})} \cong \mathbb{C}[[t_{1,1}, t_2,\ldots, t_{3g-4}]]\]
		Because the point is smooth, the map $c_\mu$ extends over it. 
		
		Consider now the $1$-dimensional smoothing family of $\pi\colon X \rightarrow\Gamma$ obtained by varying the coordinate $t_{1,1}$
		\[
		\begin{tikzcd}
			\mathcal{C}  \arrow{r}{\pi} \arrow[swap]{dr}{} & \mathcal{P} \arrow{d}{} \\
			& \Delta_{t_{1,1}}
		\end{tikzcd}
		\]
		
		Observe that $\mathcal{C}$ is smooth except for an $A_3$-singularity at either $x$ or $y$, depending on whether $k=1$ or $k=3$. By blowing-up the singularity we obtain in the central fibre of the new space $\tilde{\mathcal{C}}$ a chain $R\cup R' \cup R''$ of rational components connecting the two points $x$ and $y$. We denote by $R_2, \ldots R_{i-1}$ the rational curves of $X$ glued to $C$ at the points $x_2,\ldots, x_{i-1}$. We take the line bundle $\mathcal{O}_{\tilde{\mathcal{C}}}(Z+ \sum_{j=2}^{i-1}R_j + R)$ on $\tilde{\mathcal{C}}$ and see that its restriction to $C$ is 
		\[ L\coloneqq \mathcal{O}_C(x_2+\cdots + x_{i-1} +\frac{3-k}{2}x + \frac{k-1}{2}y - y_1-\cdots - y_i) \]
		which satisfies $L^2 \cong \mathcal{O}_C(-x-y)$. Moreover, this line bundle on $\tilde{\mathcal{C}}$ has degree 0 when restricted to a rational component of the central fibre except for the component $R'$ for which the degree is 1. In particular, by collapsing all rational components but $R'$ we obtain an element in $\Delta_0^{\mathrm{ram}}$ over the central fibre.
		
	\end{proof}
	
	We are now ready to compute the intersection of the divisor with different test curves. 
	
	\subsection{Intersection with test curves of type $A$} Before starting our computations we make a remark about some admissible covers with elliptic source curve. We do this because such covers will appear naturally in our study and the next remark will be essential in our computation. 
	\begin{remark} \label{covermultiplicity} Let $f_1\colon E \rightarrow \mathbb{P}^1$ be a map of degree $2k\geq 4$ with ramification profiles $(2k), (4,2,\ldots,2)$ and $(2,\ldots, 2)$ over three branch points. Then $\mathrm{Aut}(E,y)\neq \mathbb{Z}_6$ and $f_1 \circ j = f_1$ for any automorphism $j\colon E \rightarrow E$ fixing the point $y$ of ramification order $2k$. 
	\end{remark}
	\begin{proof} We split the problem into two cases depending on whether $k$ is even or odd. 
		
		When $k$ is even, we consider the unique map $\pi\colon \mathbb{P}^1\rightarrow \mathbb{P}^1$ with ramification profiles $(k), (2,\ldots, 2)$ and $(1,1,2,\ldots,2)$ over $\infty$, $1$ and 0, where the point of ramification order $k$ is denoted $x_1$ and the two unramified points over $0$ are denoted $x_2$ and $x_3$. Let $x_4$ be one of the other $k-1$ points of ramification order $2$. If we take the degree 2 map $g\colon E \rightarrow \mathbb{P}^1$ with branch points $x_1,x_2,x_3$ and $x_4$, the map $\pi\circ g$ has ramification profiles as in the hypothesis. Hence all maps $f_1$ are obtained in this way and it is clear that $f_1\circ j = f_1$ for the involution $j$ of the source curve $(E,y)$. 
		
		Consider the map $g\colon E \rightarrow \mathbb{P}^1$ of degree 2 with branch points $0,1,\infty$ and $x$. Then $\mathrm{Aut}(E,y) =\mathbb{Z}_4$ if and only if $x = \frac{1}{2}$. In particular, at most one of the $k-1$ maps as in the hypothesis has source curve satisfying  $\mathrm{Aut}(E,y) =\mathbb{Z}_4$. We prove that such a map exists. 
		
		We consider the map $\pi'\colon \mathbb{P}^1 \rightarrow \mathbb{P}^1$ branched over $\infty, 0$ and $1$, having ramification profiles $(\frac{k}{2}), (1,1,2,\ldots,2)$ and $(2,\ldots, 2)$ if $\frac{k}{2}$ is even or ramification profiles $(\frac{k}{2}), (1,2,\ldots,2)$ and $(1,2,\ldots, 2)$ if $\frac{k}{2}$ is odd, where we take the point of total ramification to be $\infty$ and the two unramified points to be $0$ and $1$. 
		
		Take the map $g'\colon E \rightarrow \mathbb{P}^1$ of degree $4$ having ramification profiles $(4), (4)$ and $(2,2)$ over the branch points $\infty, 0$ and $1$. Then $\pi'\circ g'\colon E \rightarrow \mathbb{P}^1$ is a map with ramification profiles as in the hypothesis and $\mathrm{Aut}(E,y) = \mathbb{Z}_4$. 
		
		Let $g\colon E \rightarrow \mathbb{P}^1$ a degree 2 map ramified over the points $0,1,\infty$ and $x$. Then $\mathrm{Aut}(E,y) = \mathbb{Z}_6$ if and only if $x$ is a primitive root of order 6. 
		
		We consider again the map $\pi\colon \mathbb{P}^1 \rightarrow \mathbb{P}^1$. We can describe it as 
		\[ \pi(t) = t(t-1)Q(t)^2\]
		where $Q(X)$ is given by the unique solution of the polynomial Pell equation 
		\[ P(X)^2 - X(X-1)Q(X)^2=1\]
		with $\mathrm{deg}(P) = deg(Q)+1 = \frac{k}{2}$. Observe that we have the initial solution $P(X) = 2X-1, Q(X) =2$ and hence our solution is given by 
		\[ P(X)-Q(X)\sqrt{X(X-1)} = (2X-1 - 2\sqrt{X(X-1)})^{\frac{k}{2}} \]
		Let $\xi$ a primitive root of order 6. Using the identity $\xi^2 -\xi = -1$ and the binomial expansion we deduce 
		\[ P(\xi) = (2\xi-1)^{\frac{k}{2}} \sum_{s=0}^{\left \lfloor{\frac{k}{4}}\right \rfloor } \frac{4^s}{3^s} \binom{\frac{k}{2}}{2s}\]
		In particular, it is clear that $\xi(\xi-1)Q(\xi)^2 = P(\xi)^2 - 1 $ is not $0$ or $1$ and hence $\xi$ is not one of the $k-1$ ramification points. 
		
		The approach when $k$ is odd is similar. We consider the unique map $\pi\colon \mathbb{P}^1\rightarrow \mathbb{P}^1$ with branch points $\infty, 0$ and $1$ having ramification profiles $(k), (1,2,\ldots, 2)$ and $(1,2,\ldots,2)$ where the point of ramification order $k$ is $\infty$ and the unramified points are $0$ and $1$ respectively. The $k-1$ maps are obtained as in the previous case. Moreover, the unicity of the map implies that if we take $\tau(t) =1-t$ the morphism fixing $\infty$ and permuting $0$ and $1$ we have 
		\[ \pi\circ\tau(t) = 1- \pi(t)\]
		implying that none of the $k-1$ other ramification points is fixed by $\tau$. Hence $\mathrm{Aut}(E,y) \neq \mathbb{Z}_4$. 
		
		To see that $\mathrm{Aut}(E,y) \neq \mathbb{Z}_6$ it is enough to show that if $\xi$ is a primitive root of order $6$, then $\pi(\xi)$ is not $0$ or $1$. The map $\pi$ can be described as 
		\[ \pi(t) = tQ(t)\]
		where $Q(X)$ is given by the unique solution of the generalized polynomial Pell equation 
		\[ XQ(X)^2 - (X-1)P(X)^2 = 1\]
		with $\mathrm{deg}(P) = \mathrm{deg}(Q) = \frac{k-1}{2}$. We see that the solution is given by 
		\[ XQ(X) - P(X)\sqrt{X(X-1)} = (X-\sqrt{X(X-1)})(2X-1-2\sqrt{X(X-1)})^{\frac{k-1}{2}}\] 
		Hence we have 
		\[ Q(X) = \sum_{s=0}^{\left \lfloor{\frac{k-1}{4}}\right \rfloor }4^sX^s(X-1)^s(2X-1)^{\frac{k-1}{2}-2s} \binom{\frac{k-1}{2}}{2s}+2(X-1)\sum_{s=0}^{\left \lfloor{\frac{k-3}{4}}\right \rfloor }4^sX^s(X-1)^s(2X-1)^{\frac{k-3}{2}-2s}\binom{\frac{k-1}{2}}{2s} \]
		As before, we get that $\xi Q(\xi)$ is not $0$ or $1$ and hence $\xi$ is not one of the $k-1$ points of ramification order $2$ of $\pi$. Hence $\mathrm{Aut}(E,y) \neq \mathbb{Z}_6$.
		
	\end{proof}
	
	Next, we compute the intersection of the divisor $\overline{D}(\mu)$ with the test curves obtained by pulling back to $\overline{\mathcal{R}}_g$ the test curve $A$. 
	\begin{proposition} \label{AAAinterA}
		We have the following intersection numbers 
		\[ \overline{D}(\mu)\cdot A_{1:g-1} = \overline{D}(\mu)\cdot A_{g-1} = 0 \textrm{ \ and \ } \overline{D}(\mu)\cdot A_1 = (3g-3)!\cdot6\cdot(2^{2i-2}- \binom{2i-1}{i})\]
	\end{proposition}
	\begin{proof}
		We first show that an element $[C\cup_{x\sim y}E, \eta] \in \overline{\mathcal{R}}_g$ with $[C,x] \in \mathcal{M}_{g-1,1}$ generic, $E$ a smooth elliptic curve and $\eta_C \not\cong \mathcal{O}_C$ is not contained in $\overline{D}(\mu)$. 
		
		We assume by contradiction that it is. Due to the genericity of $[C,x]$, the existence of a twist implies that the partition on $C$ determined by it has length $g-1$. It follows that this partition is $\mu$. 
		We know that the forgetful map 
		\[ \mathcal{H}^0_{g-1}(\mu) \rightarrow \mathcal{M}_{g-1}\]
		is finite. It follows that $x$ is not a ramification point for any of the maps $f\colon C \rightarrow \mathbb{P}^1$ with ramification profiles $\mu^{+}$ and $-\mu^{-}$ over $0$ and $\infty$. 
		
		If $[\pi\colon X \rightarrow \Gamma] \in \overline{H}_{g,\mu}$ is mapped by $c_\mu\colon\overline{H}_{g,\mu} \rightarrow \overline{\mathcal{R}}_g$ to $[C\cup_{x\sim y} E, \eta]$, the previous remarks imply $\mathrm{deg}(\pi_{|C}) = g$ and $\mathrm{ord}_x(\pi_{|C}) = 1$. This would imply $\mathrm{deg}(\pi_{|E}) = 1$, which is impossible. 
		
		Next, we describe the admissible covers $[\pi\colon X \rightarrow \Gamma]$ that are mapped by $c_\mu$ to a point of the test curve $A_1$. We denote by $x_1, \ldots, x_s$ the points of $C$ in the same fibre of $\pi_{|C}$ as $x$. As in Case II of Proposition \ref{extension} we define
		\[ c_{0,1} = \sum_{\substack{\mathbb{P} \textrm{ collapses to $x\sim y$} \\ \mathbb{P} \textrm{ is mapped to } \mathbb{P}_1}} \textrm{deg}(f_{|\mathbb{P}}\colon\mathbb{P} \rightarrow \mathbb{P}_1) \] 
		\[ c_{0,2} = \sum_{\substack{\mathbb{P} \textrm{ collapses to $x\sim y$} \\ \mathbb{P} \textrm{ is mapped to } \mathbb{P}_2}} \textrm{deg}(f_{|\mathbb{P}}\colon\mathbb{P} \rightarrow \mathbb{P}_2) \]
		and similarly $c_{j,1}$ and $c_{j,2}$ corresponding to $x_j$ for $j=\overline{1,s}$. Using the properties of admissible covers we deduce the relations 
		\[ \mathrm{ord}_x(\pi_{|C}) +c_{0,1} = \mathrm{ord}_y(\pi_{|E}) +c_{0,2}\]
		\[ \mathrm{ord}_{x_1}(\pi_{|C}) +c_{1,1} = c_{1,2} \]
		\[\ldots\]
		\[ \mathrm{ord}_{x_s}(\pi_{|C}) +c_{s,1} = c_{s,2} \] 
		
		Adding them up and using that $c_{j,1} \geq 1$ for every $j = \overline{1,s}$ we get that 
		\[\mathrm{deg}(\pi_{|C}) + s + c_{0,1} \leq 2i\]
		Together with the inequality 
		\[\mathrm{deg}(\pi_{|C}) + s  \geq 2i\]
		coming from the genericity of $[C,x]$, this implies that $c_{1,1} = \cdots = c_{s,1} = 1$, that $c_{0,1}=0$ and $\mathrm{deg}(\pi_{|C}) = 2i-s$. This implies further that $c_{1,2} = \cdots = c_{s,2} = 2$, that $\mathrm{deg}(\pi_{|E}) = 2i-2s$ and $\pi_{|E}$ has ramification profiles $(2i-2s), (4,2,\ldots,2)$ and $(2,\ldots,2)$ over $q,q_1$ and $q_2$. 
		
		Consequently, the map $\pi\colon X\rightarrow \Gamma$ is uniquely determined by the choice of a map $f\colon C \rightarrow \mathbb{P}^1$ of degree $2i-s$ with ramification order $2i-2s$ at $x$ and the choice of a map $f_1\colon E\rightarrow \mathbb{P}^1$ of degree $2i-2s$ having ramification profiles $(2i-2s), (4,2,\ldots,2)$ and $(2,\ldots,2)$ over three branch points. It follows from Theorem \ref{BrillNoether} and Proposition \ref{table} that the number of such maps $\pi\colon X \rightarrow \Gamma$ is equal to 
		\[\sum_{s=0}^{i-1}(2i-2s)\cdot \frac{(2i-1)!}{(2i-s)!s!}\cdot (i-s-1) = \frac{1}{i} \sum_{s=0}^{i-1}(i-s)(i-s-1)\binom{2i}{s}\]
		
		As a consequence of Proposition \ref{combiden}, the number of maps $\pi\colon X \rightarrow \Gamma$ with the desired properties is $2^{2i-2} - \binom{2i-1}{i}$. We observe that each choice of order for the $3g-3$ simple branch points produces a different admissible cover. We next show that each admissible cover should be counted with multiplicity 6. 
		
		Using the description of the complete local ring, we deduce that the admissible cover $\pi \colon X\rightarrow \Gamma$ admits a universal family
		\[
		\begin{tikzcd}
			\mathcal{C}  \arrow{r}{\pi} \arrow[swap]{dr}{} & \mathcal{P} \arrow{d}{} \\
			& \Delta
		\end{tikzcd}
		\] 
		where $\Delta$ is a $(3g-4)$-dimensional polydisk. This induces a map from $\Delta$ to the universal deformation of the curve $C\cup_{x\sim y}E$ and the same method as in the proof of Theorem 6 in \cite{KodMg} implies that the image intersect the singular locus $\Delta_1$ transversely at the point $[C\cup_{x\sim y}E]$. Moreover, $\mathrm{Aut}_\pi(X)$ acts on $\Delta$ and $\Delta/\mathrm{Aut}_\pi(X)$ is an open neighbourhood of $[\pi\colon X\rightarrow \Gamma]$ in $\overline{H}_{g,\mu}$. 
		
		By computing the intersection $\overline{D}(\mu)\cdot A_1$ at the level of the universal deformation of $C\cup_{x\sim y}E$, we get that each admissible cover should be counted with multiplicity $\frac{1}{2}\cdot 12 \cdot \frac{|\mathrm{Aut}(C\cup_{x\sim y}E)|}{|\mathrm{Aut}_\pi(X)|}$. Here, the factor $\frac{1}{2}$ appears because $\delta_1 = \frac{1}{2}\Delta_1$ and the factor $12$ appears as each elliptic curve shows up $12$ times in the pencil. Using Remark \ref{covermultiplicity}, we deduce that the multiplicity is always 6. We conclude that 
		\[\overline{D}(\mu)\cdot A_1 = 6\cdot(3g-3)!(2^{2i-2}- \binom{2i-1}{i})\]
		In order to conclude that the intersection is 0 with the two other test curves, we still need to show that an element $[C\cup_{x\sim y}E_{\infty}, \eta]\in \overline{\mathcal{R}}_g$ with $[E_{\infty},y]$ the singular curve in $\overline{\mathcal{M}}_{1,1}$ does not appear in the intersection. The methods we used for a smooth elliptic curve extend to this case.
		
	\end{proof}
	
	\subsection{Intersection with test curves of type $C_{g-1}$} 
	We can employ a similar approach as in the case of test curves of type $A$ for the test curves of type $C_{g-1}$. 
	\begin{proposition} We have the intersection numbers 
		\[\overline{D}(\mu)\cdot C_{1:g-1}^{g-1}= (3g-3)!\cdot(4i-4)\cdot(2i+2)\cdot\binom{2i-1}{i}, \ \overline{D}(\mu)\cdot C_{g-1}^{g-1}= (3g-3)!\cdot(4i-4)\cdot(6i-2)\cdot\binom{2i-1}{i}\]
		\[\normalfont\textrm{and \ \ } \overline{D}(\mu)\cdot C_{1}^{g-1}= (3g-3)!\cdot(4i-4)\cdot\text{\Large(}2i(4i+1)\cdot\binom{2i-1}{i} - 6(2i-1)\cdot2^{2i-2}\text{\Large)}\]	
	\end{proposition}
	\begin{proof} We start by computing the first intersection. Let $[\pi\colon X\rightarrow \Gamma]\in \overline{H}_{g,\mu}$ be an admissible cover mapped to a point $[C\cup_{x\sim y}E,\eta]$ on $C^{g-1}_{1:g-1}$. Because $c_\mu = b_\mu \circ a_\mu$ and $b_\mu$ was defined in terms of the unique twist, we distinguish three different cases for the divisorial equivalences on $C$ and $E$ implied by the twist. 
		
		\textbf{Case I:} The equivalence on $E$ is $\mathcal{O}_E(2x_1-2y) \cong \eta_E$ and the one on $C$ is \[\mathcal{O}_C(2x+x_2+\cdots+x_{i-1}-x_i-\cdots- x_{2i-1}) \cong \eta_C\] 
		There are 4 choices of $x_1$ for the first equivalence, $4\cdot (2i-1)!$ solutions for the second equivalence and one choice for each solution of a point $x$ having coefficient 2. The choices of $x_1$ are two by two identified by the involution of $E$. Moreover, as the order of the points having the same coefficient is irrelevant and as each ordering of the simple branch points produces a different admissible cover, we get $(3g-3)!\cdot2\cdot \frac{4\cdot (2i-1)!}{(i-2)!\cdot i!} \cdot 1$ elements in $\overline{H}_{g,\mu}$ having a corresponding twist as above. If $\alpha\colon X\rightarrow X$ is the automorphism acting as the involution on $E$ and fixing all other components of $X$, we get $\pi\circ \alpha \neq \pi$. 
		
		The complete local ring of $\overline{H}_{g,\mu}$ at such a point $[\pi\colon X\rightarrow \Gamma]$ is the ring of invariants of 
		\[\mathbb{C}[[t_{1,1},t_{1,2},\ldots, t_{1,i-1},t_1,\ldots, t_{3g-4}]]_{/(t_1 = t^4_{1,1} = t_{1,2}^2=\cdots = t_{1,i-1}^2)} \]
		with respect to the group $\mathrm{Aut}_\pi(X)$. This group has cardinality $2^{i-2}$ and consists only of automorphisms that fix $C$ and $E$ and act non-trivially on the $i-2$ rational components of $X$. Assuming the action to be linear, we immediately get that the complete local ring is isomorphic to 
		\[\mathbb{C}[[t_{1,1},t^2_{1,2},\ldots, t^2_{1,i-1},t_1,\ldots, t_{3g-4}]]_{/(t_1 = t^4_{1,1} = t_{1,2}^2=\cdots = t_{1,i-1}^2)} \cong \mathbb{C}[[t_{1,1},t_2,\ldots, t_{3g-4}]] \] 
		We take $\Delta$ to be the parameter space of this ring and see it is a base for a universal deformation of the map $\pi\colon X \rightarrow \Gamma$. Consequently we get a map from $\Delta$ to the universal deformation of the curve $C\cup_{x\sim y}E$. 
		
		The automorphism $\alpha\colon C\cup_{x\sim y} E \rightarrow C\cup_{x\sim y} E$ lifts to an automorphism $\alpha\colon X\rightarrow X$ identifying two by two the admissible covers. It follows that at the level of the universal deformation of the curve $C\cup_{x\sim y}E$, the branch of the image of $\Delta$ is simply tangent to the locus $\Delta_1$ parametrizing singular curves. As a consequence, all such covers appear in the count $\overline{D}(\mu)\cdot C^{g-1}_{1:g-1}$ with multiplicity $2$. A similar argument to this can be found in Lemma 3.4 in \cite{KodevenHarris1984}.

		\textbf{Case II:} The equivalence on $E$ is $\mathcal{O}_E(x_2-y) \cong \eta_E$ and the one on $C$ is 
		\[\mathcal{O}_C(2x_1+x+x_3+\cdots+x_{i-1}-x_i-\cdots- x_{2i-1}) \cong \eta_C\] 
		In this case, there is a unique choice of $x_2$ for the first equivalence, $4\cdot(2i-1)!$ solutions for the second one and $i-2$ choices for each solution of a point $x$ having coefficient $1$. The order of the points having the same coefficient is irrelevant and each ordering of the simple branch points produces a different admissible cover, hence we find $(3g-3)!\cdot \frac{4\cdot (2i-1)!}{(i-2)!\cdot i!} \cdot (i-2)$ elements in $\overline{H}_{g,\mu}$ having a corresponding twist as above. 
		
		\textbf{Case III:} The equivalence on $E$ is $\mathcal{O}_E(y-x_{2i-1})\cong \eta_E$ and the one on $C$ is
		\[ \mathcal{O}_C(2x_1+x_2+\cdots+x_{i-1}-x_i-\cdots-x_{2i-2}-x)\cong \eta_C\]
		Reasoning as in the previous cases we get $(3g-3)!\cdot\frac{4\cdot (2i-1)!}{i!(i-1)!}\cdot i$ admissible covers in $\overline{H}_{g,\mu}$ with corresponding twist as above. 
		
		Reasoning as in the proof of Theorem 6 in \cite{KodMg} we deduce that all admissible covers in the Cases II and III appear with multiplicity $1$. Hence we have 
		\[ \overline{D}(\mu)\cdot C^{g-1}_{1:g-1} = (3g-3)!\cdot \frac{4\cdot (2i-1)!}{i!\cdot(i-2)!}(2\cdot 2 + i-2+i)\]
		
		We proceed to compute the intersection $\overline{D}(\mu)\cdot C^{g-1}_{g-1}$. In this case, there are only two possibilities for the twist. 
		
		\textbf{Case I:} The equivalence on $E$ induced by the twist is $\mathcal{O}_E(2x_1-2y) \cong \eta_E$ and the one on $C$ is 
		\[ \mathcal{O}_C(2x+x_2+\cdots+x_{i-1}-x_i-\cdots - x_{2i-1}) \cong \eta_C \] 
		In this case, we have $3$ choices for $x_1$ and $4\cdot (2i-1)!$ solutions on $C$. For each such solution on $C$ there is a unique choice of the point $x$ with coefficient $2$. As the order of the points having the same coefficient is irrelevant and as each ordering of the branch points produces a different admissible cover, we get $(3g-3)!\cdot3\cdot \frac{4\cdot (2i-1)!}{i!(i-2)!}$ elements in $\overline{H}_{g,\mu}$ with corresponding twist as above. 
		
		\textbf{Case II:} The other possibility is when the equivalence on $E$ is trivial and the one on $C$ is 
		\[  \mathcal{O}_C(2x_1+x_2+\cdots+x_{i-1}-x_i-\cdots - x_{2i-1}) \cong \eta_C\]
		In this situation we get $\mathrm{deg}(\pi_{|C}) = 2i$ and hence $\mathrm{deg}(\pi_{|E}) = \mathrm{ord}_x(\pi_{|C})$. As $C$ is generic, it follows that $\mathrm{deg}(\pi_{|E}) =2$ and $x$ is one of the $6i-5$ simple ramification points of $\pi_{|C}$. The order of the points having the same coefficient in the divisorial equivalence on $C$ is irrelevant and each ordering of the $3g-3$ simple branch points produces a different admissible cover. We obtain in this way $(3g-3)!\cdot \frac{4\cdot(2i-1)!}{i!\cdot (i-2)!}\cdot(6i-5)$ admissible covers corresponding to this case. 
		
		The method in the proof of Theorem $6$ in \cite{KodMg} implies that all the admissible covers in the two cases appear with multiplicity $1$ in the intersection $\overline{D}(\mu) \cdot C^{g-1}_{g-1}$. Hence we have 
		\[\overline{D}(\mu) \cdot C^{g-1}_{g-1} = (3g-3)!\cdot \frac{4\cdot(2i-1)!}{i!\cdot(i-2)!}\cdot (3+6i-5) \] 
		
		Finally we compute $\overline{D}(\mu)\cdot C^{g-1}_1$. In this case, the twist on $C$ is trivial. Moreover, for an admissible cover $[\pi\colon X\rightarrow \Gamma]$ mapped to $C_1^{g-1}$ we have that $q_1$ and $q_2$ are contained in the component of $\Gamma$ that is the target of the elliptic curve $E$. In the notations of Proposition \ref{AAAinterA} we have the relations 
		\[ \mathrm{ord}_x(\pi_{|C}) +c_{0,1} = \mathrm{ord}_y(\pi_{|E}) +c_{0,2}\]
		\[ \mathrm{ord}_{x_1}(\pi_{|C}) +c_{1,1} = c_{1,2} \]
		\[\ldots\]
		\[ \mathrm{ord}_{x_s}(\pi_{|C}) +c_{s,1} = c_{s,2} \] 
		We distinguish two different cases depending on the position of the point over $q_1$ of ramification order $4$. 
		
		If this point is on $E$, by adding the relations we obtain that $\mathrm{deg}(\pi_{|C}) + s+ c_{0,1} \leq 2i$. The genericity of $C$ implies $\mathrm{deg}(\pi_{|C}) + s  \geq 2i-1$ and the genericity of $E$ implies $\mathrm{deg}(\pi_{|C}) + s+ c_{0,1} \leq 2i-1$. It follows that $c_{0,1} = 0$, $\mathrm{deg}(\pi_{|C}) = 2i-1-s$ and $\mathrm{ord}_x(\pi_{|C}) = 2i-2s-1$. 
		
		Consequently, the map $\pi\colon X \rightarrow \Gamma$ is uniquely determined by a map $f\colon C \rightarrow \mathbb{P}^1$ of degree $2i-s-1$ with ramification order $2i-2s-1$ at a point $x$ and a map $f_1\colon E \rightarrow \mathbb{P}^1$ having ramification profiles $(2i-2s-1,1), (4,2,\ldots,2), (2,\ldots, 2)$ and $(2,1,\ldots, 1)$ over $q,q_1,q_2$ and $q_3$, satisfying $\mathcal{O}_E(\frac{f_1^{*}(q_1)-f_1^{*}(q_2)}{2})\cong \eta_E$. 
		
		Using Theorem \ref{BrillNoether} and Proposition \ref{doiHurw}, we deduce that the number of such maps $\pi\colon X \rightarrow \Gamma$ is equal to 
		\[ \frac{1}{3} \sum_{s=0}^{i-1}4(i-s-1)(2i-2s-1)(i-s)\binom{2i-1}{s} \cdot (6i-6s-3) \]
		It follows from Proposition \ref{combiden} that this is equal to $6(i-1)(2i-1)\cdot 2^{2i-2}$ and hence the number of admissible covers $\pi\colon X\rightarrow \Gamma$ of this form is $(3g-3)!\cdot 6(i-1)(2i-1)\cdot 2^{2i-2}$. If we consider the automorphism $\alpha\colon X\rightarrow X$ acting as the involution on $E$ and fixing all other components, we see that $\pi\circ \alpha \neq \pi$. Reasoning as in Case I of the computation $\overline{D}(\mu)\cdot C^{g-1}_{1:g-1}$, we deduce that all the admissible covers above should be counted with multiplicity $2$. 
		
		The other possible case is when the point of ramification order $4$ is on a rational component collapsing to $x_j$ when we stabilize $X$. In this case we have $c_{j,1}\geq 2$. This implies the inequality $\mathrm{deg}(\pi_{|C}) + s + c_{0,1} \leq 2i-1$. Using this and the inequality $\mathrm{deg}(\pi_{|C}) + s \geq 2i-1$ coming from the genericity of $C$ we get $c_{0,1} = 0$ and $\mathrm{deg}(\pi_{|C}) =2i-1-s$. In this case the ramification orders of $\pi_{|C}$ at $x$ and $x_j$ are $2i-2s-2$ and $2$ respectively. 
		
		Consequently, the map $\pi\colon X \rightarrow \mathbb{P}^1$ is uniquely determined by a degree $2i-s-1$ map $f\colon C\rightarrow \mathbb{P}^1$ having ramification orders $2i-2s-2$ and $2$ at two points $x$ and $x_j$ in the same fibre and a map $f_1 \colon E \rightarrow \mathbb{P}^1$ having ramification profiles $(2i-2s-2), (2,\ldots, 2), (2,\ldots,2)$ and $(2,1,\ldots, 1)$ over four branch points $q, q_1, q_2, q_3$ and satisfying $\mathcal{O}_E(\frac{f_1^{*}(q_1)-f_1^{*}(q_2)}{2})\cong \eta_E$. 
		
		From Theorem \ref{BrillNoether}, the number of such maps $f\colon C \rightarrow \mathbb{P}^1$ is $8 s(i-s-1)(i-s)(4i-4s-5) \cdot \binom{2i-1}{s}$. From Proposition \ref{ellipticsix}, the number of such maps $f_1\colon E \rightarrow \mathbb{P}^1$ is $2$. It follows that the number of such maps $\pi\colon X \rightarrow \Gamma$ is equal to 
		\[ 16 \sum_{s=0}^{i-1}s(i-s-1)(i-s)(4i-4s-5)\cdot\binom{2i-1}{s}\]
		which using Proposition \ref{combiden} we compute to be 
		\[ 8i(i-1)(4i+1)\cdot\binom{2i-1}{i} - 36(2i-1)(i-1)\cdot 2^{2i-2} \] 
		Again, every ordering of the simple branch points produces a different point in $\overline{H}_{g,\mu}$ and hence the number of admissible covers in this case is the number we computed multiplied by $(3g-3)!$. We show that all such admissible covers should be counted with multiplicity $1$. 
		
		The complete local ring of $\overline{H}_{g,\mu}$ at the point $[\pi\colon X \rightarrow \Gamma]$ is the ring of invariants of 
		\[ \mathbb{C}[[t_{1,1},t_{1,2},t_1,\ldots, t_{3g-4}]]_{/(t_1 = t_{1,1}^{2i-2s-2} = t_{1,2}^2)} \] 
		The method of Remark \ref{covermultiplicity} can be employed to prove $\pi \circ \alpha = \pi$ where $\alpha\colon X \rightarrow X$ is the automorphism acting as the involution on $E$ and as identity on all other components. The group $\mathrm{Aut}_\pi(X)$ is generated by $\alpha$ and the automorphism $\tau$ acting non-trivially on the components collapsing to the point of order $2$ and fixing the other components. 
		
		By considering the ring of invariants with respect to $\tau$ we get that the complete local ring of $\overline{H}_{g,\mu}$ at the point $[\pi\colon X \rightarrow \Gamma]$ is the ring of invariants of 
		\[ \mathbb{C}[[t_{1,1},t_{1,2}^2, t_1,\ldots, t_{3g-4}]]_{/(t_1 = t_{1,1}^{2i-2s-2} = t_{1,2}^2)} \cong \mathbb{C}[[t_{1,1}, t_1,\ldots, t_{3g-4}]]_{/(t_1 = t_{1,1}^{2i-2s-2})} \] 
		with respect to $\alpha$. The same method of mapping the parameter space to the universal deformation of $C\cup_{x\sim y} E$ and proceeding as in \cite{KodMg} implies that all admissible covers appear with multiplicity $1$. It follows that 
		\[\normalfont \overline{D}(\mu)\cdot C_{1}^{g-1}= \cdot(3g-3)!\cdot\text{\Large(}8i(i-1)(4i+1)\cdot\binom{2i-1}{i} - 24(2i-1)(i-1)\cdot2^{2i-2}\text{\Large)}\]
	\end{proof}
	
	\subsection{Intersection with test curves of type $B$} Finally we compute the intersection of our divisor with the test curves $B'$ and $B''$. The work of understanding which admissible covers $[\pi\colon X\rightarrow \Gamma]$ map to these test curves already appears implicitly in Proposition \ref{extension}. We are left with the task of computing their number and their multiplicities. 
	\begin{proposition} We have the following intersection numbers: 
		\[\overline{D}(\mu)\cdot B' = (3g-3)!\cdot 8(i^2-i)\cdot (2^{4i-2}-1)\cdot\binom{2i-1}{i} \textrm{ \ \ and } \]
		\[\normalfont \overline{D}(\mu)\cdot B'' = (3g-3)!\text{\Large [}(8i^3-8i^2-2i)\binom{2i-1}{i} - (2i-1)(6i-8)2^{2i-2}\text{\Large ]}\]
	\end{proposition}
	\begin{proof}
		We start by computing $\overline{D}(\mu)\cdot B'$. The admissible covers in $\overline{H}_{g,\mu}$ mapping to $\Delta_0'$ are described in Case I and Case IV of Proposition \ref{extension}. 
		
		Let $[\pi\colon X\rightarrow \Gamma]$ be an admissible cover as in Proposition \ref{extension}, Case I, mapped to $B'$. The number of solutions for the divisorial equivalence 
		\[ \mathcal{O}_C(2x_1+x_2+\cdots+x_{i-1}-x_i-\cdots- x_{2i-1})\cong \eta_C \]
		is equal to $4\cdot (2i-1)!$ for any of the $2^{4i-2}-1$ elements $\eta_C \in \mathrm{Pic}(C)[2]\setminus \left\{0\right\} $. It follows there are $\frac{(2i-1)!}{i!(i-2)!}(2^{4i}-4)$ choices of a map $\pi_{|C}$ having degree $2i$ and ramification profiles $(4,2,\ldots,2)$ and $(2,\ldots,2)$ over two points $q_1$ and $q_2$. For the generic point $x$ on $C$, there are $2i-1$ choices of a point $y$ in the same fiber of $\pi_{|C}$. Moreover, there is a rational component $R$ of $X$ passing through $x$ and $y$ and mapping $2:1$ to $\mathbb{P}_2$. 
		
		We fix a map $\pi_{|C}$ and a point $y$ as just discussed. Two orderings of the $(3g-3)$ simple branch points produce the same admissible cover if and only if they differ by transposing the order of the two branch points on $\mathbb{P}_2$. Hence we get 
		\[ \frac{(3g-3)!}{2}\cdot \frac{(2i-1)!}{i!(i-2)!}\cdot (2i-1)\cdot (2^{4i}-4)\]
		distinct admissible covers. It is immediate from the description in the proof of Proposition \ref{extension} that all these covers are counted with multiplicity $2$. 
		
		For an admissible cover $[\pi\colon X \rightarrow \Gamma]$ as in case IV of Proposition \ref{extension}, the map $\pi_{|C}$ has degree $2i$ and ramification profiles $(2,\ldots,2)$ and $(2,\ldots,2)$ over $q$ and $q_2$, with the generic point $x$ one of the ramified points over $q$. Moreover, there is a rational component $R$ mapping $4:1$ to $\mathbb{P}_2$ connecting $x$ with one of the other $i-1$ points in the same fiber over $q$ as $x$. The component $R$ contains the point of ramification order $4$ over $q_1$. 
		
		Every ordering of the $(3g-3)$ simple branch points produces a different admissible cover. We obtain in this way 
		\[ (3g-3)!\cdot \frac{(2i-1)!}{i!(i-1)!}\cdot (i-1)\cdot (2^{4i-2}-1)\]
		admissible covers and we deduce from Proposition \ref{extension}, Case IV that all should be counted with multiplicity $4$. 
		
		It follows that 
		\[\overline{D}(\mu)\cdot B' = (3g-3)!\cdot 8(i^2-i)\cdot (2^{4i-2}-1)\cdot\binom{2i-1}{i}  \]
		
		Next we compute $\overline{D}(\mu)\cdot B''$. In Case II of Proposition \ref{extension} we outlined three possible types of admissible covers $[\pi\colon X \rightarrow \Gamma]$ mapping to $B''$. 
		
		The first type is when $\pi_{|C}$ has degree $2i-s-1$ and ramification order at $x$ and $y$ equal to $i-s$ and $i-s-1$. In this case, for the rational component $R$ joining $x$ and $y$, the map $\pi_{|R}$ has ramification profiles $(i-s,i-s-1), (4,2,2,\ldots,2)$ and $(2,\ldots,2)$ over $q,q_1$ and $q_2$. 
		
		From Theorem \ref{BrillNoether}, the number of choices of such a map $\pi_{|C}$ is \[(i-s-1)\text{\large[}(i-s-1)(2i-2s-1)-1\text{\large]}\cdot \binom{2i-1}{s}\] 
		Proposition \ref{table} implies the choice of $\pi_{|R}$ is unique. Moreover, each ordering of the simple branch points produces a different admissible cover and the discussion in Proposition \ref{extension} implies each of them should be counted with multiplicity $2i-2s-1$. 
		
		The contribution to the count coming from this case is 
		\[\normalfont (3g-3)!\cdot \sum_{s=0}^{i-1} (2i-2s-1)(i-s-1)\text{\Large[}(i-s-1)(2i-2s-1)-1\text{\Large]}\cdot \binom{2i-1}{s} \]
		which we deduce from Proposition \ref{combiden} to be equal to 
		\[  (3g-3)!\cdot \text{\Large[} \frac{3}{2} \cdot (2i-1)(i-1)\cdot 2^{2i-2} - 2(i-1)i\cdot \binom{2i-1}{i}\text{\Large]} \] 
		
		The second type is when $\pi_{|C}$ has degree $2i-s-1$ and ramification orders $x$ and $y$ equal to $i-s-1$ and $i-s$. The number of such maps is 
		\[(i-s)\text{\large[}(i-s)(2i-2s-1)-1\text{\large]}\cdot \binom{2i-1}{s}\]
		Other than that, everything follows identically as in the previous case and we get a contribution of 
		\[(3g-3)!\cdot \sum_{s=0}^{i-1} (2i-2s-1)(i-s)\text{\large[}(i-s)(2i-2s-1)-1\text{\large]}\cdot \binom{2i-1}{s}  \]
		which we deduce from Proposition \ref{combiden} to be equal to 
		\[  (3g-3)!\cdot \text{\Large[} \frac{3}{2} \cdot (2i-1)(i-1)\cdot 2^{2i-2} + 2(i-1)i\cdot \binom{2i-1}{i}\text{\Large]} \] 
		
		The third type is when $\pi_{|C}$ is a map of degree $2i-s-1$ having ramification orders $i-s-1, i-s-1$ and $2$ at $x,y$ and another point in the same fiber of $\pi_{|C}$ as $x$ and $y$. For the rational component $R$ joining $x$ and $y$, we have that $\pi_{|R}$ has ramification profiles $(i-s-1, i-s-1), (2,\ldots,2)$ and $(2,\ldots,2)$ over $q, q_1$ and $q_2$. 
		
		The number of maps $\pi_{|C}$ of this type is equal to 
		\[8s(i-s)(i-s-1)^2\binom{2i-1}{s} - 2s(i-s-1)(i-s+1)\binom{2i-1}{s}\]
		
		Each ordering of the simple branch points produces a different admissible cover and each appears in $\overline{D}(\mu)\cdot B''$ with multiplicity $2$. Hence the contribution in this case is 
		\[ (3g-3)!\sum_{s=0}^{i-1}16s(i-s)(i-s-1)^2\binom{2i-1}{s} - (3g-3)!\sum_{s=0}^{i-1}4s(i-s-1)(i-s+1)\binom{2i-1}{s}\]
		From Proposition \ref{combiden} we deduce the identities 
		\[\sum_{s=0}^{i-1}16s(i-s)(i-s-1)^2\binom{2i-1}{s} = 8(i-1)i^2\binom{2i-1}{i}-8(i-1)(2i-1)\cdot 2^{2i-2}\]
		\[ \sum_{s=0}^{i-1}4s(i-s-1)(i-s+1)\binom{2i-1}{s} = 2i\cdot\binom{2i-1}{i} + (2i-1)(i-3)\cdot 2^{2i-2} \]
		
		Putting everything together, we conclude that 
		\[\normalfont \overline{D}(\mu)\cdot B'' = (3g-3)!\text{\Large [}(8i^3-8i^2-2i)\binom{2i-1}{i} - (2i-1)(6i-8)2^{2i-2}\text{\Large ]}\]
	\end{proof}
	
	\subsection{Conclusions (Proof of Theorem \ref{maintheorem}):} We denote 
	\[ \overline{D}(\mu)  = (3g-3)!\cdot \text{\Large(} a\lambda - b_0'\delta_0' - b_0''\delta_0'' - b_0^{\mathrm{ram}}\delta_0^{\mathrm{ram}}- b_1\delta_1 - b_{g-1}\delta_{g-1} - b_{1:g-1}\delta_{1:g-1}- \cdots \text{\Large)} \]
	From the test curve computations we deduce we have the following system of equations 
	\[ a - 4b_0'-4b_0^{\mathrm{ram}}+b_{1:g-1} = a -12b_0'+b_{g-1}= 0, \ \ a-4b_0''-4b_0^{\mathrm{ram}} +b_1 = 2\cdot 2^{2i-2} - 2\cdot \binom{2i-1}{i}\]
	\[b_{1:g-1} = (2i+2)\cdot\binom{2i-1}{i}, \ \ b_{g-1} = (6i-2)\binom{2i-1}{i}, \ \ b_1 = 2i(4i+1)\cdot\binom{2i-1}{i}-6(2i-1)\cdot 2^{2i-2}\]
	\[ (8i-4)b_0' - b_{g-1}-b_{1:g-1}=(8i^2-8i)\cdot\binom{2i-1}{i} \ \ \mathrm{and} \] 
	\[(4i-2)b_0''-b_1 = (8i^3-8i^2-2i)\binom{2i-1}{i} - (2i-1)(6i-8)\cdot2^{2i-2}\]
	
	This is a solvable system of 8 equations in 7 unknowns. We compute the coefficients to be 
	\[ a = \frac{12i^2+10i-2}{2i-1}\cdot\binom{2i-1}{i}, \ \ b_0' = \frac{2i^2}{2i-1}\cdot\binom{2i-1}{i}\]
	\[b_0^{\mathrm{ram}} = \frac{2i^2+3i-1}{2i-1}\cdot\binom{2i-1}{i} \ \ \mathrm{and} \ \ b_0'' = \frac{4i^3}{2i-1}\cdot\binom{2i-1}{i}-(3i-1)\cdot 2^{2i-2}\]
	
	We remark that in Theorem \ref{maintheorem}, the contribution coming from the order of the $3g-3$ simple branch points is not taken into account. \hfill $\square$
	\subsection{A divisor in $\overline{\mathcal{R}}_{2i+1}$}
	For genus $g= 2i+1$ and partition $\mu = (2,\ldots,2,-2,\ldots,-2)$ of length $g-1$, we can apply the same procedure to compute the divisor $\overline{D}(\mu)$. This is the divisor $\overline{D}_{2i+1:2}$ appearing in \cite{FarLud}. By the method of test curves we deduce this divisor has the following coefficients: 
	\[a = \frac{1}{2i-1}\binom{2i}{i}\cdot (3i+1), \ \ b_0' = \frac{1}{2i-1}\binom{2i}{i}\cdot \frac{i}{2},  \] 
	\[b_0'' = \frac{1}{2i-1}\binom{2i}{i}\cdot i^2, \ \ b_0^{\mathrm{ram}} =\frac{1}{2i-1}\binom{2i}{i}\cdot \frac{2i+1}{4}\] 
	
	As a consequence, we get the intersection $\overline{D}(\mu) \cdot B^{\mathrm{ram}}_0 = (3g-3)!\cdot\binom{2i}{i}\cdot(2^{2g-3}-2)$ and by describing explicitly the points and their multiplicity as in the proof of Proposition \ref{extension} we deduce 
	\begin{corollary}
		The degree of the map  
		\[ \mathcal{H}^0_{2i}(\underbrace{2,\ldots, 2}_{i \ \mathrm{entries}}, \underbrace{-2,\ldots, -2}_{i-1 \ \mathrm{entries}}, -1, -1) \rightarrow \mathcal{M}_{2i,1}\]
		forgetting all but the last marking is equal to $(2i)!\cdot(2^{4i-2}-1)$. 
	\end{corollary}
	
	We observe that the coefficient $b_0''$ of the divisor $\overline{D}(\mu)$ differs from the one computed in \cite{FarLud}. This happens because the map 
	\[\phi: \wedge^i \mathcal{H} \otimes \mathcal{A}_{0,0} \rightarrow \mathcal{A}_{i-1,1}\] 
	used in \cite{FarLud} to compute this divisor degenerates above the locus $\Delta_0''$. We recall that fiberwise, $\phi$ is given over a point $[X,\eta]$ as 
	\[ \wedge^iH^0(X, \omega_X) \otimes H^0(X,\omega_X\otimes\eta)\rightarrow H^0(X, \wedge^{i-1}M_{X}\otimes \omega_X^2\otimes \eta)\] 
	where $M_X$ is the Lazarsfeld vector bundle of $\omega_X$ and $M_L$ denotes the Lazarsfeld vector bundle of the line bundle $L$, as in \cite{Laza}.
	
	If $[C_{/x\sim y}, \eta]$ is a generic point of $\Delta_0''$ it follows that
	\[ \wedge^iH^0(C,\omega_C(x+y)) \otimes H^0(C,\omega_C(x+y))\rightarrow H^0(X, \wedge^{i-1}M_{\omega_C+x+y}\otimes \omega_C^2(2x+2y))\]  
	is not an isomorphism. We use the exact sequence 
	\[ 0 \rightarrow \wedge^iM_{\omega_C+x+y} \otimes \omega_C(x+y) \rightarrow \wedge^iH^0(C,\omega_C(x+y))\otimes\omega_C(x+y) \rightarrow \wedge^{i-1}M_{\omega_C+x+y}\otimes \omega_C^2(2x+2y) \rightarrow 0\]
	to deduce that $h^0(C, \wedge^iM_{\omega_C+x+y} \otimes \omega_C(x+y))$ and $h^1(C, \wedge^iM_{\omega_C+x+y} \otimes \omega_C(x+y))$ are not 0. 
	
	Using Proposition 1.3.3 in \cite{Laza} regarding the Green-Lazarsfeld property $(N_{i-1})$ we deduce
	\begin{proposition} \label{greenlazarsfeld}
		Let $g = 2i$ and $[C,x,y]$ a generic element in $\mathcal{M}_{g,2}$. Then $\omega_C(x+y)$ fails to satisfy the property $(N_{i-1})$.  
	\end{proposition}
	
	While deriving this result using test curves is an interesting approach of Proposition \ref{greenlazarsfeld}, this result is not new. It immediately follows from \cite{AproduVoisin} Th\'eor\`eme 0.3 and \cite{Farsyzy} Theorem 3.7. 
	
	
	%
	%

	\textbf{Acknowledgements} I am grateful to Gavril Farkas and Johannes Schmitt for their insightful remarks and ideas, that proved to be fundamental in obtaining the results of this article. I am also thankful to Carlos Maestro P\'erez and Marian Aprodu for their remarks, which helped me simplify and better explain some of the arguments. I also want to thank the anonymous referee for the comments and suggestions.  \ \\
	\textbf{Data availability statement}: not applicable.

	%
	%

\bibliography{main}
\bibliographystyle{alpha}
\end{document}